\documentclass[onefignum,onetabnum]{siamart250211}

\usepackage{appendix}

\usepackage{enumerate}
\usepackage{multirow}
\usepackage{subfig}
\usepackage{tabularx}
\usepackage{makecell}
\usepackage{booktabs}
\usepackage{graphicx}
\usepackage{subcaption}
\usepackage{float}
\usepackage{comment}

\usepackage[normalem]{ulem}
\usepackage{cancel}
\newcommand{\GPLadd}[1]{\textcolor{black}{#1}}

\newcommand{\Toadd}[1]{\textcolor{black}{#1}}
\usepackage{bm}
\usepackage{amsmath}
\usepackage{amstext}
\usepackage{amssymb}
\usepackage{bbm}
\usepackage{amsbsy}
\usepackage{mathtools}

\newsiamremark{example}{Example}
\newsiamremark{remark}{Remark}

    \newcommand{\beq}{\begin{equation}}
    \newcommand{\eeq}{\end{equation}}
    \newcommand{\beqnn}{\begin{equation*}}
    \newcommand{\eeqnn}{\end{equation*}}
        
    \def\bb {\boldsymbol{b}}

    \def\bp {\boldsymbol{p}}
    
    \def\bu {\boldsymbol{u}}

    \def\bw {\boldsymbol{w}}
    \def\bx {\boldsymbol{x}}
    \def\by {\boldsymbol{y}}
    \def\bz {\boldsymbol{z}}

    \newcommand{\matr}[1]{\bm{#1}}

    \newcommand{\norm}[1]{\left\Vert{#1} \right\Vert}

    \def\R {\mathbb{R}}
    \def\Rn {\R^n}

    \def\prox {\text{prox}}
\providecommand{\freg}{f_{\mathrm{reg}}}
    
    \newcommand{\diff}{\mathop{}\!d}
    \DeclareMathOperator*{\argmin}{arg\,min}

    \newcommand{\normsq}[1]{\left\Vert{#1} \right\Vert_{2}^{2}}
    
    \newcommand{\normone}[1]{\left\Vert{#1} \right\Vert_{1}}
    \newcommand{\norminf}[1]{\left\Vert{#1} \right\Vert_{\infty}}

    \def\JBVS  {J_{\text{BVS}}}
    \def\JPAM {J_{\text{PAM}}}
    \def\JPQM {J_{\text{PQM}}}
    \def\JNN   {J_{\text{NN}}}
    \def\psiNN {\psi_{\text{NN}}}
    \def\PsiMP {\Psi_{\text{MP}}}
\providecommand{\bx}{\boldsymbol{x}}
\providecommand{\by}{\boldsymbol{y}}
\providecommand{\bz}{\boldsymbol{z}}
\providecommand{\bw}{\boldsymbol{w}}
\providecommand{\bu}{\boldsymbol{u}}

\providecommand{\freg}{f_{\mathrm{reg}}}
\providecommand{\psiNN}{\psi_{\mathrm{NN}}}
\providecommand{\JNN}{J_{\mathrm{NN}}}
\providecommand{\JBVS}{J_{\mathrm{BVS}}}
\providecommand{\prox}{\operatorname{prox}}
\providecommand{\uPM}{u_{\mathrm{PM}}}
\providecommand{\TV}{\mathrm{TV}}
\providecommand{\Seps}{S_{\epsilon}}
\providecommand{\Keps}{K_{\epsilon}}
\providecommand{\coloneqq}{\mathrel{:=}}
\providecommand{\Dstar}{\mathsf{D}_{\star}}
\providecommand{\Dpir}{\mathsf{D}_{\mathrm{PIRATE}}}
\providecommand{\Dplus}{\mathsf{D}_{\mathrm{PIRATE+}}}

\title{Deep learning methods for inverse problems using connections between proximal operators and Hamilton--Jacobi equations}

\author{Oluwatosin Akande\thanks{Industrial and Systems Engineering,
Lehigh University,
200 West Packer Avenue,
Bethlehem, PA 18015, USA, (\url{oaa323@lehigh.edu})} \and Gabriel P. Langlois\thanks{Department of Mathematics, University of Illinois Urbana-Champaign, Chicago, IL, USA (\url{gp42@illinois.edu}).} \and Akwum Onwunta\thanks{Industrial and Systems Engineering,
Lehigh University,
200 West Packer Avenue,
Bethlehem, PA 18015, USA, (\url{ako221@lehigh.edu}).}}

\begin{document}
\maketitle

\begin{abstract}
Inverse problems are important mathematical problems that seek to recover model parameters from noisy data. Since inverse problems are often ill-posed, they require regularization or incorporation of prior information about the underlying model or unknown variables. \GPLadd{Proximal operators, ubiquitous in nonsmooth optimization, are central to this because they encode priors and yield efficient iterative algorithms.} They have also recently become key to modern machine learning methods, e.g., plug-and-play methods with learned denoisers and deep neural architectures for learning priors of proximal operators. The latter was developed partly due to recent work characterizing proximal operators of nonconvex priors as subdifferentials of convex potentials. In this work, we propose to leverage connections between proximal operators and Hamilton--Jacobi partial differential equations (HJ PDEs) to develop deep learning architectures for learning the prior. In contrast to other existing methods, we learn the prior directly without recourse to inverting the prior after training. We present numerical results in dimensions up to $64$, where the recovered prior is evaluated in a single forward pass.
\end{abstract}

\section{Introduction}\label{sec:intro}
Inverse problems seek to recover model parameters from noisy data. Applications in science and engineering include image recovery from noisy measurements, deblurring, tomographic reconstruction, and compressive sensing~\cite{AEOV2023,arridge2019solving, bertero2021introduction,isakov2017inverse}. Since inverse problems are often ill-posed, it is essential to include regularization or prior information about the underlying model or unknown variables. Proximal operators are central to this: they provide a flexible and computationally convenient way to encode priors and to build efficient iterative algorithms (e.g., proximal (sub)gradients and splitting methods). More recently, proximal operators have become key ingredients for state-of-the-art machine learning methods, such as plug-and-play methods that replace explicit regularizers by learned denoisers~\cite{hu2023plug,jia2025plug} or deep neural architectures parameterizing proximal maps or their gradients\GPLadd{~\cite{fang2024whats}}.  

Formally, the proximal operator of a proper function $J \colon \Rn \to \R\cup\{+\infty\}$ is defined via observed data $\bx \in \Rn$, a parameter $t>0$, and the minimization problem
\begin{equation}\label{eq:def_prox_min}
	S(\bx,t) = \min_{\by \in \Rn} \left\{\frac{1}{2t}\normsq{\bx-\by} + J(\by)\right\}.
\end{equation}
\GPLadd{The value function defined by~\cref{eq:def_prox_min} is the Moreau envelope of $tJ$~\cite[Definition~12.20]{bauschke2017convex}.} The proximal operator is the set-valued function
\begin{equation}\label{eq:def_prox_argmin}
\prox_{tJ}(\bx) = \argmin_{\by \in \Rn} \left\{\frac{1}{2t}\normsq{\bx-\by} + J(\by)\right\}.
\end{equation}
Here, $t$ controls the trade-off between the quadratic data-fidelity term and the prior. \GPLadd{Borrowing language from statistics, we call $J$ a \emph{prior} of the operator $\prox_{tJ}$.} In practice one often works directly with $\prox_{tJ}$ rather than the prior. 

The recent work of Gribonval and Nikolova~\cite{gribonval2020characterization} in nonsmooth optimization has extended the characterization of proximal operators with convex priors to those with nonconvex priors, showing that they are subdifferentials of certain convex potentials. These properties, in particular, were used in~\cite{fang2024whats} to develop learned proximal networks (LPNs), deep learning methods that learn the prior of a proximal operator from data.

The paper~\cite{gribonval2020characterization} did not, however, discuss the well-established, existing connections between proximal operators and Hamilton--Jacobi Partial Differential Equations (HJ PDEs)~\cite{chaudhari2018deep,darbon2015convex,darbon2021bayesian,darbon2016algorithms,osher2023hamilton}. To see these connections, consider the following HJ PDE, whose Hamiltonian is quadratic and whose initial data is the prior $J$:
\begin{equation} \label{eqn:intro2}
\begin{dcases} 
\frac{\partial S}{\partial t}(\bx,t)+\frac{1}{2}\normsq{\nabla_{\bx}S(\bx,t)}= 0, &\ \GPLadd{(\bx,t)} \in \Rn \times(0,+\infty),\\
S(\bx,0)=J(\bx), &\ \bx \in \Rn.
\end{dcases}
\end{equation}
If $J$ is uniformly Lipschitz continuous, then the unique \emph{viscosity solution} of the HJ PDE is given by~\cref{eq:def_prox_min}. Moreover, at a point of differentiability $\bx$, there holds
\begin{equation}\label{eqn:early_prox}
\prox_{tJ}(\bx) = \bx - t\nabla_{\bx}S(\bx,t).
\end{equation}
Finally, the viscosity solution is \emph{semiconcave}, that \GPLadd{is}, the function $\bx \mapsto \frac{1}{2}\normsq{\bx} - tS(\bx,t)$ is convex. Pairing this with~\cref{eqn:early_prox}, we see that $\prox_{tJ}$ is obtained by differentiating a convex function. This connects proximal operators to HJ PDEs, a connection already established in the literature, and to the stronger\footnote{To the best \GPLadd{of }our knowledge, this characterization result was unknown in the theory of HJ PDEs.} characterization obtained in~\cite{gribonval2020characterization}.

In this paper, we propose to leverage the theory of viscosity solutions of HJ PDEs to develop deep learning methods for learning the prior function $J$ in~\cref{eq:def_prox_argmin} directly from data. To describe our approach, suppose the solution $(\bx,t) \mapsto S(\bx,t)$ to the HJ PDE~\cref{eqn:intro2} is known. We consider later the case in which only samples of the solution are available. This problem was investigated in~\cite{barron1999regularity,claudel2011convex, colombo2020initial, esteve2020inverse, misztela2020initial}. In particular,~\cite{esteve2020inverse} showed that when $\bx \mapsto S(\bx,t)$ is uniformly Lipschitz continuous and $\bx \mapsto \frac{1}{2}\normsq{\bx} - tS(\bx,t)$ is convex, there exists a prior $J$ recovering $S(\bx,t)$ exactly. Moreover, there is a natural choice for the prior, obtained by reversing the time in the HJ PDE~\cref{eqn:intro2} and using $(\bx,t) \mapsto S(\bx,t)$ as the terminal condition. The resulting \emph{backward viscosity solution} yields the prior $\JBVS\colon \Rn \to \R$ which admits the representation formula
\begin{equation}\label{eq:prior_backward}
\JBVS(\by) = \sup_{\bx \in \Rn} \left\{S(\bx,t) - \frac{1}{2t}\normsq{\bx-\by}\right\}.
\end{equation}
Here, $J(\by) \geqslant \JBVS(\by)$ for every $\by \in \Rn$, with $\JBVS(\by) = J(\by)$ whenever $\by = \bx - t\nabla_{\bx}S(\bx,t)$, where $\bx$ is a point of differentiability of $\bx \mapsto S(\bx,t)$. Moreover, 
\[
\inf_{\by \in \Rn} \left\{\frac{1}{2t}\normsq{\bx-\by} + \JBVS(\by)\right\} = S(\bx,t) \ \text{for every $\bx \in \Rn$.}
\]
Thus the prior $\JBVS$ recovers the function $\GPLadd{\bx} \mapsto S(\bx,t)$, although in general $\prox_{tJ}$ and $\prox_{t\JBVS}$ may not agree everywhere. Nonetheless, this provides a principled way to estimate the prior, at least when $S(\bx,t)$ is known.

We shall focus in this paper on the case where $\bx \mapsto S(\bx,t)$ is unknown but \GPLadd{we }have access to some samples $\{\bx_{k},S(\bx_{k},t),\nabla_{\bx}S(\bx_{k},t)\}_{k=1}^{K}$ (or some subset thereof) with $t$ fixed. We propose to learn the prior $\by \mapsto \JBVS(\by)$ by leveraging the crucial fact that $\by \mapsto \GPLadd{t}\JBVS(\by) + \frac{1}{2}\normsq{\by}$ is convex, thus enabling approaches based on deep learning and convex neural networks.

\GPLadd{\noindent
\textbf{Contributions of this paper:} \textbf{(i)} We formulate the recovery of the prior underlying a proximal operator from data as an inverse problem of recovering the backward viscosity solution $\JBVS$ of an appropriate set of Hamilton--Jacobi equations. After clarifying connections between proximal maps of nonconvex priors and the theory of backward viscosity solutions of Hamilton--Jacobi PDEs, we propose to learn the backward viscosity solution directly via convex neural networks by leveraging its semiconvexity property. 
\textbf{(ii)} We prove a transfer bound,~\cref{prop:transfer_values}, under which a potential fitted to within $\varepsilon_{\psi}$ on a region $B$ yields by conjugation a prior within $\varepsilon_{\psi}/t$ of the backward viscosity solution $\JBVS$, provided $B$ contains the two points at which the conjugates are attained. The training region must therefore cover the preimages of the query points. The piecewise quadratic minorant of $t\JBVS$ built from $K$ well-spread samples errs by $\Theta(K^{-2/n})$ (\cref{cor:jbvs_approx} and~\cref{eq:pqm_fill}), and the minima of $k$ quadratics escape that rate (\cref{ex:quad_exact}).
\textbf{(iii)} We test the framework in dimensions $2$ through $64$ on four priors known in closed form, namely $\normone{\bx}$, a minimum of two quadratics, $-\normsq{\bx}/4$ and $-\normone{\bx}$, scoring each recovery by its relative $L^{2}$ error against $\JBVS$. Our recovery evaluates the prior in one forward pass, while the inversion approach of~\cite{fang2024whats} optimizes at every query point. It stays within a factor of $3.2$ of that inversion through $d = 8$ and beats it on all four priors at $d = 64$. On the two priors of~\cref{ex:quad_exact}, conjugating convex max-affine quadratic potentials (defined in~\cref{subsec:maxplus}) to recover the prior directly outperforms the inversion approach by up to five orders of magnitude. We then recover the regularizer implicit in a total-variation posterior-mean denoiser. Following the work of~\cite{darbon2021bayesian}, who showed that a total variation posterior mean denoiser is the proximal mapping of some implicit prior, we show numerically in~\cref{subsec:tv_exp} that the implicit prior can be successfully learned from data. Finally, we perform numerical experiments on the two Plug-and-Play operators released in~\cite{hu2023plug} in~\cref{subsec:pnp_exp}.}

\textbf{Related work:}  Hamilton--Jacobi PDEs are important to many scientific and engineering applications such as optimal control~\cite{Bardi1997Optimal,fleming2006controlled,mceneaney2006max,parkinson2018optimal} and physics~\cite{Caratheodory1965CalculusI,Caratheodory1967CalculusII}, inverse problems for imaging sciences~\cite{darbon2015convex,darbon2021bayesian,darbon2019decomposition,darbon2022hamilton}, optimal transport~\cite{meng2024primal,onken2021ot}, game theory~\cite{BARRON1984213,Evans1984Differential,ruthotto2020machine}, and machine learning~\cite{chen2024leveraging,zou2024leveraging}. Recent works focus on developing specialized methods for solving high-dimensional HJ PDEs via, e.g., representation formulas or deep learning. These specialized methods leverage certain properties of HJ PDEs, including stochastic aspects and representation formulas~\cite{bardi1998hopf,darbon2022hamilton,darbon2016algorithms,mceneaney2006max}, to approximate solutions to HJ PDEs more accurately and efficiently than general-purpose methods. See, e.g.,~\cite{darbon2023neural,darbon2020overcoming,darbon2021some,meng2022sympocnet,park2025neural} for recent works along these lines and~\cite{meng2025recent} for a review of state-of-the-art numerical methods for HJ PDEs.

Deep learning methods have become popular for computing solutions to high-dimensional PDEs and their inverse problems. This is because neural networks can be trained on data to approximate high-dimensional, nonlinear functions using efficient optimization algorithms. They have been used to approximate solutions to PDEs without discretizing on a numerical grid, which is one reason they can mitigate the curse of dimensionality. There is a fairly comprehensive literature on deep learning methods for solving PDEs in general, e.g., see~\cite{beck2020overview,cuomo2022scientific,karniadakis2021physics}.

\textbf{Organization of this paper:} We present background information on proximal operators, Hamilton--Jacobi equations, and convex neural networks in~\cref{sec:prelim}. Next, we discuss recent results concerning the inverse problem for Hamilton--Jacobi equations when the solution is available, and how they relate to proximal operators and learning priors in inverse problems, in~\cref{sec:backward}. We present some theoretical results in~\cref{sec:inverse_incomplete_info}, where we study the inverse problem for Hamilton--Jacobi equations when only incomplete information is available about its solution. 
We suggest via arguments from max-plus algebra theory for Hamilton--Jacobi PDEs how to learn from data the solution to a certain \GPLadd{Hamilton--Jacobi} terminal value problem, which can then be used as an estimate for learning the prior function in a proximal operator. 
We present in~\cref{sec:numerics} some numerical experiments for learning the initial data of certain Hamilton--Jacobi PDEs using convex neural networks and the theory of inverse Hamilton--Jacobi PDEs. Finally, we summarize our results in~\cref{sec:discussion}.
\section{Background}\label{sec:prelim}
We present here some background on proximal operators, HJ PDEs, connections between them, and convex neural networks. For comprehensive references, we refer the reader to~\cite{cannarsa2004semiconcave,evans2022partial,rockafellar2009variational}.

\subsection{Proximal operators}\label{subsec:background_convex}
Let $J\colon \Rn \to \R\cup\{+\infty\}$ denote a proper function (i.e., $J(\bx) < +\infty$ for some $\bx \in \Rn$ and $J(\bx) > -\infty$ for every $\bx \in \Rn$). Consider the minimization problem $(\bx,t) \mapsto S(\bx,t)$ defined in~\cref{eq:def_prox_min} and its proximal operator $(\bx,t) \mapsto \prox_{tJ}(\bx)$ defined in~\cref{eq:def_prox_argmin}. We say a proper function $f_{t}\colon \Rn \to \GPLadd{\Rn}$ is a \emph{proximal operator} of $tJ$ if $f_{t}(\bx) \in \prox_{tJ}(\bx)$ for every $\bx \in \Rn$. Gribonval and Nikolova~\cite{gribonval2020characterization} proved that proximal operators are characterized in terms of the function $\psi\colon \Rn \times [0,+\infty) \to \R\cup\{+\infty\}$ defined by
\begin{equation}\label{eq:psi_defn}
    \psi(\bx,t) = \frac{1}{2}\normsq{\bx} - tS(\bx,t).
\end{equation}
\begin{theorem}\label{thm:char_prox}
A proper function $f_{t}\colon \Rn\to \Rn$ is a proximal operator of $tJ$ if and only if $\bx \mapsto \psi(\bx,t)$ is proper, lower semicontinuous, and convex and $f_{t}(\bx) \in \partial_{\bx}\psi(\bx,t)$. Moreover, $f_{t}$ is uniformly Lipschitz continuous with constant $L>0$ if and only if $\bx \mapsto (1 - 1/L)\normsq{\bx}/2 + tJ(\bx)$ is proper, lower semicontinuous and convex.
\end{theorem}
\begin{proof}
	See~\cite[Theorem 3 and Proposition 2]{gribonval2020characterization}
\end{proof}
The characterization of proximal operators in~\cref{thm:char_prox} is closely related to the concepts of \emph{semiconcave} and \emph{semiconvex} functions.
\begin{definition}\label{defn:semiconcavity}
Let $\mathcal{C} \subset \Rn$. We say $g \colon \mathcal{C} \to \R$ is $C$-semiconcave with $C\geqslant 0$ if it is continuous and
    \[
    \lambda g(\bx_{1}) + (1-\lambda)g(\bx_{2}) - g(\lambda\bx_{1} + (1-\lambda)\bx_{2}) \leqslant \lambda(1-\lambda)\GPLadd{\frac{C}{2}}\normsq{\bx_{1}-\bx_{2}}
    \]
    for every $\bx_{1},\bx_{2} \in \mathcal{C}$ such that $\lambda\bx_{1} + (1-\lambda)\bx_{2} \GPLadd{\in} \mathcal{C}$ and $\lambda \in [0,1]$. We say $g$ is semiconvex if $-g$ is semiconcave.
\end{definition}
\begin{remark}\label{rem:semiconcavity}
It can be shown~\cite[Chapter 1]{cannarsa2004semiconcave} that a function $g$ is $C$-semiconcave with $C\geqslant 0$ if and only if $\bx \mapsto g(\bx) - \frac{C}{2}\normsq{\bx}$ is concave, if and only if $g = g_{1} + g_{2}$, where $g_{1}$ is concave and $g_{2} \in C^{2}(\Rn)$ with $\norminf{\nabla_{\bx}^{2}g_{2}} \leqslant C$.
\end{remark}
Combining formula~\cref{eq:psi_defn},~\cref{defn:semiconcavity} and~\cref{rem:semiconcavity}, we find $\bx \mapsto \psi(\bx,t)$ is convex if and only if $\bx \mapsto tS(\bx,t)$ is semiconcave. We will see later that semiconcavity is an important concept in the theory of HJ PDEs for characterizing their generalized solutions. But before moving on to present some background on HJ PDEs, we give below an instructive example.
\begin{example}[The negative absolute value prior]\label{ex:negl1_1}
Let $J(x) = -|x|$ and consider the one-dimensional problem
\[
S(x,t) = \min_{y \in \R} \left\{\frac{1}{2t}(x-y)^2 - |y|\right\}.
\]
A global minimum $y^{*}$ of this problem satisfies the first-order optimality condition
\[
0 \in (y^{*}-x)/t - \partial|y^{*}| \iff y^{*} \in \begin{cases}
    x + t, &\ \text{if $y^{*}>0$,} \\ 
    [x-t, x+t] &\ \text{if $y^{*}=0$},\\
    x - t, &\ \text{if $y^{*}<0$.}
\end{cases}
\]
A short calculation shows that
\begin{equation}\label{eq:ex_hj_sol}
S(x,t) = -\frac{t}{2} - |x| \quad \text{and} \quad \prox_{tJ}(x) = \begin{cases}
    x + t, &\ \text{if $x>0$,} \\ 
    \{-t, t\} &\ \text{if $x=0$},\\
    x - t, &\ \text{if $x<0$.}
\end{cases}
\end{equation}
Thus, a selection $f_{t}(x) \in \prox_{tJ}(x)$ differs only at $x = 0$. In any case, the function $x \mapsto \psi(x,t)$ in \cref{thm:char_prox} and its subdifferential $x \mapsto \partial_{x}\psi(x,t)$ are given by 
\[
\psi(x,t) = \frac{1}{2}x^2 - tS(x,t) = \frac{1}{2}x^{2} + t|x| + \frac{t^2}{2} \quad \text{and} \quad \partial_{x}\psi(x) = \begin{cases}
x + t,\, &\text{if $x>0$},\\
[-t,t],\, &\text{if $x=0$},\\
x - t,\, &\text{if $x<0$}.
\end{cases}
\]
We see that any selection $f_{t}(x) \in \prox_{tJ}(x)$ satisfies $f(x) \in \partial\psi(x,t)$.
\end{example}

\subsection{Hamilton--Jacobi Equations}\label{subsec:background_hjpde}
In this section, we briefly review some elements of the theory of HJ PDEs, including the method of characteristics, viscosity solutions of HJ PDEs, and the Lax--Oleinik formula, and discuss how these concepts connect to proximal operators.
The discussion is not comprehensive; see~\cite{evans2022partial} and references therein for a more detailed treatment.

To ease the presentation, we consider only the first-order HJ PDE~\cref{eqn:intro2}.
\subsubsection{Characteristic equations}
The characteristic equations of~\cref{eqn:intro2} are given by the dynamical system
	\begin{equation}\label{eq:char_general}
	\begin{cases}
        \dot{\bx}(t) &= \bp(t),\\
		\dot{\bp}(t) &= 0, \\
		\dot{\bz}(t) &= \frac{1}{2}\normsq{\bp(t)},
	\end{cases}
	\end{equation}
where $\bz(t) = S(\bx(t),t)$ and $\GPLadd{\bz(0)} = J(\bx(0))$. Here, $t \mapsto \bp(t)$ is constant with $\bp(t) \equiv \bp(0) \in \Rn$. The characteristic line that arises from $\bx(0) \in \Rn$ is $\bx(t) = \bx(0) +t\bp(0)$, and so $\bz(t) = \bz(0) \GPLadd{{}+\frac{t}{2}}\normsq{\bp(0)}$.  Combining these facts, we find
\[
    S(\bx(t),t) = J(\bx(0)) + \GPLadd{\frac{t}{2}}\normsq{\bp(0)}.
\]
If we write $\bx(t) \equiv \bx$ and $\bp(0) = \nabla_{\bx}S(\bx,t)$, assuming formally that the spatial gradient exists at $\bx$, then $\bx(0) = \bx - t\nabla_{\bx}S(\bx,t)$, and we find the representation
\begin{equation}\label{eq:rep_formula_char}
    S(\bx,t) = \GPLadd{\frac{t}{2}}\normsq{\nabla_{\bx}S(\bx,t)} + J(\bx - t\nabla_{\bx}S(\bx,t)).
\end{equation}
This gives an implicit relation among $S$, its spatial gradient, and the initial data $J$. Next, we turn to the explicit representation of solutions to~\cref{eqn:intro2}.

\subsubsection{Viscosity solutions and the Lax--Oleinik formula}\label{subsubsec:viscosity_solns}
The initial value problem~\cref{eqn:intro2} (and HJ PDEs with general Hamiltonians) may not have a unique generalized solution, that is, a function satisfying the HJ PDE almost everywhere together with the initial condition $S(\bx,0) = J(\bx)$. 
\begin{example}
Let $J \equiv 0$ in~\cref{eqn:intro2} and take $n=1$. The corresponding HJ PDE has infinitely many solutions: For instance, the functions $S_{1}$ and $S_{2}$ given by
\[
S_{1}(x,t) = 0, \quad S_{2}(x,t) = \begin{cases}
0,\, &\text{if $|x| \geqslant \GPLadd{t/2}$}, \\
x-\GPLadd{t/2},\, &\text{if $0\leqslant x \leqslant \GPLadd{t/2}$},\\
-x-\GPLadd{t/2},\, &\text{if $\GPLadd{-t/2} \leqslant x \leqslant 0$},
\end{cases}
\]
satisfy $S_{1}(x,0) = S_{2}(x,0) = 0$ and solve the corresponding HJ PDE \GPLadd{almost everywhere. Indeed, on the set where $0 < |x| < t/2$ we have $\partial_{t}S_{2} = -\tfrac{1}{2}$ and $\tfrac{1}{2}|\partial_{x}S_{2}|^{2} = \tfrac{1}{2}$, so $S_{2}$ satisfies the PDE off the measure-zero set $\{x=0\} \cup \{|x| = t/2\}$. More generally, for each $a \in (0,1]$, the function $S_{a}$ defined by $S_{a}(x,t) = a|x| - a^{2}t/2$ if $|x| \leqslant at/2$ and by $S_{a}(x,t) = 0$ otherwise is a generalized solution with the same initial data}.
\end{example}

The notion of \emph{viscosity solution} was introduced in~\cite{crandall1983viscosity} to solve this problem. Under appropriate conditions (see~\cite{bardi1998hopf,crandall1992user,crandall1983viscosity}), the viscosity solution is unique and admits a representation formula. Specifically, for the initial value problem~\cref{eqn:intro2} with uniformly Lipschitz continuous initial data $J$, the unique viscosity solution is given by the Lax--Oleinik formula (with quadratic Hamiltonian)
\begin{equation}\label{eq:lax_formula}
	S(\bx,t) = \inf_{\by \in \Rn} \left\{\frac{1}{2t}\normsq{\bx-\by} + J(\by)\right\}.
\end{equation}
The unique viscosity solution has two important properties. First, the function $\bx \mapsto S(\bx,t)$ is (1/t)-semiconcave, that is, the function $\bx \mapsto \psi(\bx,t)$ defined in~\cref{eq:psi_defn} is convex, as per~\cref{thm:char_prox}. Second, at any point of differentiability of $\bx \mapsto S(\bx,t)$,
\begin{equation}\label{eq:minimum_representation}
\nabla_{\bx}S(\bx,t) = \frac{\bx-f_{t}(\bx)}{t} \iff f_{t}(\bx) = \bx - t\nabla_{\bx}S(\bx,t),
\end{equation}
where $f_{t}(\bx)$ denotes a global minimizer in~\cref{eq:lax_formula}. Note that substituting this expression into formula~\cref{eq:rep_formula_char} obtained from the characteristic equations yields~\cref{eq:lax_formula}, as expected. \GPLadd{For convex $J$, formula~\cref{eq:minimum_representation} is the standard expression for the gradient of the Moreau envelope~\cite[Proposition~12.30]{bauschke2017convex}; we give the derivation from the characteristic equations only to keep the later sections self-contained.}
\begin{example}[The negative absolute value prior, continued.]
    Let $J(x) = -|x|$ in the HJ PDE~\cref{eqn:intro2}. The function $J$ is uniformly Lipschitz continuous and, as such, the Lax--Oleinik formula $S(x,t) = -\frac{t}{2} - |x|$ is the unique viscosity solution of the corresponding HJ PDE. Note $\bx \mapsto S(\bx,t)$ is differentiable everywhere except at $\bx = 0$ and $\prox_{tJ}(\bx) = \bx - t\nabla_{\bx}S(\bx,t)$ everywhere except at $\bx = 0$ (see~\eqref{eq:ex_hj_sol}).
\end{example}

In summary, \GPLadd{whenever} the function $(\bx,t) \mapsto S(\bx,t)$ is the viscosity solution of the HJ initial value problem~\cref{eqn:intro2}\GPLadd{, a proper function $f_{t}$ is a proximal operator of $tJ$ precisely when $f_{t}(\bx)$ is a global minimizer in the Lax--Oleinik formula~\cref{eq:lax_formula} at every $\bx \in \Rn$}. The minimization problem underlying $\prox_{tJ}(\bx)$ is exactly the Lax--Oleinik representation formula of the viscosity solution of~\cref{eqn:intro2}. We will see in the next section how to leverage these connections for learning the prior when $\bx \mapsto J(\bx)$ is not available but $(\bx,t) \mapsto S(\bx,t)$ is available. But before proceeding, we briefly review convex neural networks, which will be used later in this work.

\subsection{Convex neural networks}\label{subsec:background_nns}
Convex neural networks, specifically input convex neural networks (ICNNs), were introduced by Amos et al.~\cite{amos2017input} to allow for the efficient optimization of neural networks within structured prediction and reinforcement learning tasks. The core premise of an ICNN is to constrain the network architecture such that the output is a convex function with respect to the input. To achieve convexity, the network typically employs a recursive structure for $k = 0,\dots,j-1$
\begin{equation}
    \bz_{k+1} = g(\matr{W}_k \bz_k + \matr{H}_k \by + \bb_k), \quad f(\by;\theta) = \bz_j,
\end{equation}
where $\by$ and $\bz_k$ are the input to the network and the hidden features at layer $k$, the matrices $\matr{W}_k$ and $\matr{H}_k$ and the vector $\bb_k$ are the parameters of layer $k$, and $g$ is the activation function. To guarantee the convexity of the output with respect to the input $\by$, specific constraints are imposed on the parameters and the activation function, which are (i) the weights $\matr{W}_k$, which connect the previous hidden layer to the current one, must be non-negative ($\matr{W}_k \geqslant 0$), and (ii) the activation function $g$ must be convex and nondecreasing~\cite{fang2024whats}.

Fang et al.~\cite[Proposition~3.1]{fang2024whats} use the ICNN architecture and the characterization of proximal operators to develop learned proximal networks (LPNs) for inverse problems. LPNs require stricter conditions than standard ICNNs. While standard ICNNs often use ReLU activation, LPNs require the activation function $g$ to be twice continuously differentiable. This smoothness is essential to ensure that the proximal operator is the gradient of a twice continuously differentiable function~\cite[Theorem 2]{gribonval2020characterization}. Consequently, LPNs typically utilize smooth activations like the softplus function, a $\beta$-smooth approximation of ReLU~\cite[Section 3]{fang2024whats}.


\section{Connections between learning priors and the inverse problem for Hamilton--Jacobi equations}\label{sec:backward}
In this section, we discuss the inverse problem of learning the prior in the proximal operator~\cref{eq:def_prox_argmin}: given $t>0$ and some function $\bx \mapsto S(\bx,t)$, assess whether there exists a prior function $J$ that can recover $\bx \mapsto S(\bx,t)$ and, if so, estimate it. Due to the connections between proximal operators and HJ equations, as discussed in Subsections~\ref{subsec:background_convex}--\ref{subsec:background_hjpde}, our starting point will be to discuss the inverse problem from the point of view of HJ equations.

We summarize in the next subsection some of the main results for this problem, based on the results of~\cite{esteve2020inverse} and other related works~\cite{claudel2011convex,colombo2020initial,misztela2020initial}.

\subsection{Reachability and inverse problems for Hamilton--Jacobi equations}\label{subsec:reachability}
We consider the inverse problem associated to the HJ initial value problem~\cref{eqn:intro2}: given $t>0$ and a function $(\bx,t) \mapsto S(\bx,t)$, find the set of initial data $J \colon \Rn \to \R$ such that the viscosity solution of~\cref{eqn:intro2} is $S(\bx,t)$. That is, we wish to find the set
\begin{equation}\label{eq:reachable_init_conds}
\begin{alignedat}{1}
I_{t}(S) &\coloneqq \{\text{$J \colon \Rn \to \R$ is uniformly Lipschitz continuous} \\
&\qquad\qquad : \text{$S(\bx,t)$ is obtained from~\eqref{eqn:intro2} at time $t$}\}.
\end{alignedat}
\end{equation}
We say the function $(\bx,t) \mapsto S(\bx,t)$ is \emph{reachable} if the set $I_{t}(S)$ is nonempty. The main reachability result for the initial value problem~\cref{eqn:intro2} is the following:
\begin{theorem}\label{thm:reachability}
    Suppose $\bx \mapsto S(\bx,t)$ is uniformly Lipschitz continuous. Then the set $I_{t}(S)$ defined in~\cref{eq:reachable_init_conds} is nonempty if and only if $\bx \mapsto tS(\bx,t)$ is \GPLadd{$1$-}semiconcave\GPLadd{, equivalently if and only if the function $\bx \mapsto \psi(\bx,t)$ defined in~\cref{eq:psi_defn} is convex}.
\end{theorem}
\begin{proof}
    This follows from~\cite[Theorem 2.2, Theorem 6.1, and Definition 6.2]{esteve2020inverse}.
\end{proof}

Now, assume $(\bx,t) \mapsto S(\bx,t)$ is reachable. \GPLadd{We now describe the set $I_{t}(S)$.} Since we evolve forward in time $J$ from $0$ to $t$ via~\cref{eqn:intro2} to obtain $S(\bx,t)$, a natural approach then is to do the opposite: evolve backward in time $S(\bx,t)$ from $t$ to $0$ via an appropriate HJ PDE. That is, we consider the terminal value problem 
\begin{equation} \label{eq:finite_time_backward}
\begin{dcases} 
\frac{\partial \bw}{\partial \tau}(\by,\tau) + \frac{1}{2}\normsq{\nabla_{\by}\bw(\by,\tau)}= 0 &\ (\by,\tau) \in \Rn \times [0,t),\\
\bw(\by,t)=S(\by,t), &\ \by \in \Rn.
\end{dcases}
\end{equation}
Under certain conditions, the HJ PDE~\cref{eq:finite_time_backward} has a unique viscosity solution:
\begin{theorem}\label{thm:prior_reconstruction}
  Suppose $\bx \mapsto S(\bx,t)$ is uniformly Lipschitz continuous and \GPLadd{$(1/t)$-}semiconcave. Then the HJ PDE~\cref{eq:finite_time_backward} has a unique viscosity solution given by
  \begin{equation}\label{eq:backward_hj_prior}
        \bw(\by,\tau) = \sup_{\bx \in \Rn} \left\{S(\bx,t) - \GPLadd{\frac{1}{2(t-\tau)}}\normsq{\bx-\by}\right\}.
  \end{equation}
Moreover, the function $\by \mapsto \GPLadd{(t-\tau)}\bw(\by,\tau)$ is semiconvex with unit constant.
\end{theorem}
\begin{proof}
See~\cite[Section 4, Equation 4.4.2]{barron1999regularity} and~\cite[Chapter 1]{cannarsa2004semiconcave}.
\end{proof}

The viscosity solution of~\eqref{eq:finite_time_backward} is sometimes called the \emph{backward viscosity solution} (BVS) to distinguish it from its forward counterpart given by~\cref{eqn:intro2}. The BVS at $\tau = 0$ corresponds to fully evolving backward in time the function $\bx \mapsto S(\bx,t)$. In what follows, we write $\JBVS \coloneqq \bw(\cdot,0)$. We can use~\cref{eq:psi_defn} to write
\begin{equation}\label{eq:formula_recovery_prior}
t\JBVS(\by) + \frac{1}{2}\normsq{\by} = \sup_{\bx \in \Rn} \left\{\langle\bx,\by\rangle - \psi(\bx,t)\right\}.
\end{equation}
The \GPLadd{right-hand} side is the convex conjugate $\psi^{*}(\cdot,t)$ evaluated at $\GPLadd{\by}$, which is well-defined because $\bx \mapsto \psi(\bx,t)$ is proper, lower semicontinuous and convex.

Theorem~\ref{thm:prior_reconstruction} suggests that $\JBVS$ is an initial condition that can reach $\bx \mapsto S(\bx,t)$. The next result stipulates that this is correct and that it is ``optimal'', in the sense that it bounds from below any other reachable initial condition $J \in I_{t}(S)$.
\begin{theorem}\label{thm:prior_optimality}
    Let $\JBVS$ denote the solution of the backward HJ terminal value problem~\cref{eq:finite_time_backward} at time $\tau = 0$. Then $J \in I_{t}(S)$ if and only if
    \[
    J(\by) \geqslant \JBVS(\by) \ \text{for every $\by \in \Rn$, with equality for every $\by \in X_{t}(S)$, where}
    \]
    \[
    X_{t}(S) \coloneqq \left\{\bx - t\nabla_{\bx}S(\bx,t) : \text{$\bx \mapsto S(\bx,t)$ is differentiable at $\bx \in \Rn$}\right\}.
    \]
\end{theorem}
\begin{proof}
    See~\cite[Theorems 2.3 and 2.4]{esteve2020inverse}.
\end{proof}
\Cref{thm:prior_optimality} is a consequence of the semiconcavity of $\bx \mapsto S(\bx,t)$, which regularizes the backward viscosity solution of~\cref{eq:finite_time_backward}. We illustrate this below with the negative absolute value prior.
\begin{example}[The negative absolute value prior, continued.]\label{ex:negl1_2}
    Let $J(x) = -|x|$ in the HJ PDE~\cref{eqn:intro2}. Recall that its unique viscosity solution is given by the Lax--Oleinik formula $S(x,t) = -\frac{t}{2} - |x|$. We now compute the corresponding unique backward viscosity solution to the terminal-value problem~\cref{eq:finite_time_backward}. The solution is well-defined because $\bx \mapsto S(\bx,t)$ is uniformly Lipschitz continuous and concave. We have
    \[
    \JBVS(x) = \sup_{\by\in\R} \left\{-\frac{t}{2} - |y| - \frac{1}{2t}(x-y)^2\right\} = -\frac{t}{2} - \inf_{\by \in \R} \left\{\frac{1}{2t}(x-y)^2 + |y|\right\}.
    \]
    The infimum on the \GPLadd{right-hand} side corresponds to the proximal operator of the function $\by \mapsto |y|$, which is the soft-thresholding operator:
    \begin{equation}\label{eq:soft_thresholding}
    \argmin_{y \in \R} \left\{\frac{1}{2t}(x-y)^2 + |y|\right\} = \begin{cases}
    x - t,\, &\text{if $x > t$}, \\
    0,\, &\text{if $x \in [-t,t]$}, \\
    x + t,\, &\text{if $x < -t$}.
\end{cases}
    \end{equation}
    This gives
    \[
    \JBVS(x) = \begin{cases}
        -|x|,\, &\text{if $x<-t$ or $x>t$}, \\
        -\frac{t}{2} - \frac{x^2}{2t},\, &\text{if $x \in [-t,t]$}.
    \end{cases}
    \]
    Here, a simple calculation shows $\mathcal{X}_{t}(S) = (-\infty,-t] \cup [t, +\infty)$, and we find $J(x) > \JBVS(x)$ on $(-t,t)$, as expected from~\cref{thm:prior_optimality}. Moreover,
    \begin{equation}\label{eq:true_jbvs}
    t\JBVS(x) + \frac{1}{2}x^2 = \begin{cases}
        \frac{1}{2}(x-t)^2 - \frac{t^2}{2},\, &\text{if $x>t$}, \\
        -\frac{t^2}{2},\, &\text{if $x \in [-t,t]$}, \\
        \frac{1}{2}(x+t)^2 - \frac{t^2}{2},\, &\text{if $x<-t$},
    \end{cases}
    \end{equation}
    and we observe $x \mapsto t\JBVS(x) + \frac{1}{2}x^2$ is convex, as expected from~\cref{thm:prior_reconstruction}.
\end{example}
The results here apply when the function $\bx \mapsto S(\bx,t)$ is known. \GPLadd{In the next section, we consider the setting in which only a finite set of values of $\bx \mapsto S(\bx,t)$ is available.}


\section{Learning priors and the inverse problem for Hamilton--Jacobi Equations with incomplete data}\label{sec:inverse_incomplete_info}
We consider here the inverse problem associated to the HJ initial value problem~\cref{eqn:intro2} with incomplete data: given $t>0$ and samples $\{\bx_{k},S(\bx_{k},t),\nabla_{\bx}S(\bx_{k},t)\}_{k=1}^{K}$, learn the prior $J$ that ``best'' recovers the observations $\{S(\bx_{k},t),\nabla_{\bx}S(\bx_{k},t)\}_{k=1}^{K}$. Assuming $\bx \mapsto S(\bx,t)$ is uniformly Lipschitz continuous, recall from~\cref{thm:reachability} that $\bx \mapsto S(\bx,t)$ is reachable if and only if it is semiconcave. If reachable, the backward viscosity solution $\by \mapsto \JBVS(\by)$ of the HJ terminal value problem~\cref{eq:finite_time_backward} recovers $(\bx,t) \mapsto S(\bx,t)$. Thus, we will focus on learning $\JBVS$.

We start with two observations. First, \GPLadd{let $\bx_{k}$ be a sample at which $\bx \mapsto S(\bx,t)$ is differentiable. Since $\bx \mapsto S(\bx,t)$ is uniformly Lipschitz continuous, Rademacher's theorem places almost every sample in this set, and differentiability at $\bx_{k}$ makes the minimizer in the Lax--Oleinik formula~\cref{eq:lax_formula} unique. That minimizer is given by}~\cref{eq:minimum_representation}:
\begin{equation}\label{eq:yk_formula}
S(\bx_{k},t) = \frac{1}{2t}\normsq{\bx_{k}-\by_{k}} + J(\by_{k}),\ \text{where} \ \by_{k} = \bx_{k} - t\nabla_{\bx}S(\bx_{k},t).
\end{equation}
Moreover, \cref{thm:prior_optimality} implies $J(\by)\geqslant\JBVS(\by)$ for every $\by \in \Rn$, \GPLadd{with equality on the set $X_{t}(S)$. The point $\by_{k}$ of~\cref{eq:yk_formula} lies in $X_{t}(S)$ by construction, so} we have $J(\by_{k}) = \JBVS(\by_{k})$. Taken together, this gives samples $\{\by_{k},J(\by_{k})\}_{k=1}^{K}$ of the true underlying prior $J$ (and $\JBVS$). Second, recall the function $\by \mapsto \GPLadd{t}\JBVS(\by) + \frac{1}{2}\normsq{\by}$ is convex by \cref{thm:prior_reconstruction}. Thus instead of learning $\JBVS$, we can equivalently learn the convex function $\by \mapsto \GPLadd{t}\JBVS(\by) + \frac{1}{2}\normsq{\by}$ using the samples $\{\by_{k},J(\by_{k})\}_{k=1}^{K}$.

Before we propose methods for learning this convex function numerically using deep learning and convex neural networks (the topic of~\cref{sec:numerics}), we discuss below different approaches for learning $\JBVS$, in particular by approximating $\by \mapsto \GPLadd{t}\JBVS(\by) + \frac{1}{2}\normsq{\by}$ from below as a piecewise convex function at the points $\left\{\by_{k}\right\}_{k=1}^{K}$. We consider this below in~\cref{subsec:piecewise_approxs}. This approximation problem is related to max-plus algebra theory for approximating solutions to HJ PDEs~\cite{akian2006max,fleming2000max,gaubert2011curse}; we discuss this in~\cref{subsec:maxplus} and specialize these results to the problem of approximating $\by \mapsto \GPLadd{t}\JBVS(\by) + \frac{1}{2}\normsq{\by}$ in general in~\cref{subsec:maxplus_inverse_hjpdes}.

\subsection{Piecewise approximations}\label{subsec:piecewise_approxs}
We consider here piecewise approximations of the prior $\by \mapsto \JBVS(\by)$ using the samples $\{\bx_{k},S(\bx_{k},t),\nabla_{\bx}S(\bx_{k},t)\}_{k=1}^{K}$ and formula~\cref{eq:yk_formula}. We consider first using a piecewise affine minorant (PAM) approximation, and then, assuming some regularity on $\JBVS$, using a piecewise quadratic minorant (PQM) approximation.

\subsubsection{Piecewise affine approximation}
We first consider the PAM approximation of the convex function $\by \mapsto t\JBVS(\by) + \frac{1}{2}\normsq{\by}$:
\begin{equation}\label{eq:PAMS}
t\JPAM(\by) + \frac{1}{2}\normsq{\by} \coloneqq \max_{k \in \{1,\dots,K\}} \left\{t\JBVS(\by_{k}) + \frac{1}{2}\normsq{\by_{k}} + \left\langle\bx_{k}, \by - \by_{k}\right\rangle\right\}.
\end{equation}
Then $\JPAM(\by) \leqslant \JBVS(\by)$ for every $\by \in \Rn$, with $\JPAM(\by_{k}) = \JBVS(\by_{k})$ at each $k \in \{1\,\dots,K\}$. A short calculation gives
\begin{equation}\label{eq:PAMS_explicit}
t\JPAM(\by) = \max_{k \in \{1,\dots,K\}} \left\{t\JBVS(\by_{k}) + \frac{1}{2}\normsq{\bx_{k} - \by_{k}} - \frac{1}{2}\normsq{\bx_{k}-\by}\right\}.
\end{equation}

\GPLadd{We now evaluate $\JPAM$ as an initial condition for~\cref{eqn:intro2}.} A short calculation yields
\begin{equation*}
\begin{alignedat}{1}
\inf_{\by \in \Rn} &\left\{\frac{1}{2t}\normsq{\bx-\by} + \JPAM(\by)\right\} = \\
&\qquad\qquad\qquad\begin{cases}
    \GPLadd{\frac{1}{t}\left(\frac{1}{2}\normsq{\bx} + h_{K}(\bx)\right)} &\ \text{if $\GPLadd{\bx \in \operatorname{conv}\{\bx_{1},\dots,\bx_{K}\}}$}, \\
    \GPLadd{-\infty}, &\ \text{otherwise},
\end{cases}
\end{alignedat}
\end{equation*}
\GPLadd{where
\[
h_{K}(\bx) \coloneqq \max\left\{\sum_{k=1}^{K}\lambda_{k}\left(tS(\bx_{k},t) - \tfrac{1}{2}\normsq{\bx_{k}}\right) \;:\; \lambda \in \Delta_{K},\ \sum_{k=1}^{K}\lambda_{k}\bx_{k} = \bx\right\}
\]
and $\Delta_{K} \coloneqq \{\lambda \in \R^{K} : \lambda_{k} \geqslant 0 \text{ for every } k,\ \sum_{k=1}^{K}\lambda_{k} = 1\}$. Since the function $\bx \mapsto tS(\bx,t) - \frac{1}{2}\normsq{\bx} = -\psi(\bx,t)$ is concave, we have $h_{K}(\bx_{k}) = tS(\bx_{k},t) - \frac{1}{2}\normsq{\bx_{k}}$ at every sample.} \GPLadd{The PAM approximation thus recovers $\{S(\bx_{k},t)\}_{k=1}^{K}$, but away from the samples it returns only a piecewise linear interpolation of the concave function $-\psi(\cdot,t)$, which diverges to $-\infty$ outside the convex hull of the samples.} It is therefore generally inadequate to use as a meaningful learned prior or estimate of $\JBVS$.

\subsubsection{Piecewise quadratic approximation}
Assume now $\by \mapsto t\JBVS(\by)$ is semiconvex with constant $(1-\alpha)$, $\alpha >0$\GPLadd{, and fix a margin $\alpha' \in (0,\alpha)$}. Then $\by \mapsto t\JBVS(\by) + \frac{1}{2}\normsq{\by}$ is $\GPLadd{\alpha}$\GPLadd{-}strongly convex, and we can approximate it via its piecewise quadratic minorant (PQM) \GPLadd{of margin $\alpha'$}:
    \[
	\begin{alignedat}{1}
	t\JPQM(\by) + \frac{1}{2}\normsq{\by} &\coloneqq \max_{k \in \{1,\dots,K\}} \biggl\{t\JBVS(\by_{k}) \ + \\
    &\qquad\qquad \left. \frac{1}{2}\normsq{\by_{k}} + \langle\bx_{k},\by-\by_{k}\rangle +\frac{\GPLadd{\alpha'}}{2}\normsq{\by-\by_{k}} \right\}.
	\end{alignedat}
    \]
	Then, $\JPQM(\by) \leqslant \JBVS(\by)$ for every $\by \in \Rn$, with $\JPQM(\GPLadd{\by_{k}}) = \JBVS(\by_{k})$ at each $k \in \{1,\dots,K\}$. Moreover, a short calculation gives
    \begin{equation*}
	t\JPQM(\by) =  \max_{k \in \{1,\dots,K\}} \left\{\GPLadd{tJ(\by_{k})} + \frac{1}{2}\normsq{\bx_{k}-\by_{k}} - \frac{1}{2}\normsq{\bx_{k}-\by} + \frac{\GPLadd{\alpha'}}{2}\normsq{\by-\by_{k}}\right\}.
	\end{equation*}

\GPLadd{We evaluate $\JPQM$ in the same way.} A formal calculation, carried out in~\cref{app:pqm_calc}, yields
\begin{equation}\label{eq:sol_hj_pde_pqm}
    \inf_{\by \in \Rn} \left\{\frac{1}{2t}\normsq{\bx-\by} + \JPQM(\by)\right\} = \GPLadd{J(\by_{k}) + {}}\frac{1}{2t}\normsq{\bx-\by_{k}} \GPLadd{{}-\frac{1}{2t\GPLadd{\alpha'}}\normsq{\bx-\bx_{k}}}
\end{equation}
for some $k \in \{1,\dots,K\}$\GPLadd{; in particular, at $\bx = \bx_{k}$ the right-hand side reduces to $J(\by_{k}) + \frac{1}{2t}\normsq{\bx_{k}-\by_{k}} = S(\bx_{k},t)$ by~\cref{eq:yk_formula}}. The PQM approximation now leads to an approximation of $(\bx,t) \mapsto S(\bx,t)$ that recovers $\{S(\bx_{k},t)\}$ and is finite everywhere. In the next section, we describe how max-plus algebra theory~\cite{akian2006max,fleming2000max,gaubert2011curse} can be used to quantify the approximation errors more precisely.

\subsection{Max-plus algebra theory for Hamilton--Jacobi PDEs and approximation results}\label{subsec:maxplus}
We consider here max-plus algebra techniques for approximating solutions to certain HJ PDEs. Let $\alpha > 0$ and let $\Psi\colon \Rn \to \R$ denote a $(1-\alpha)$-semiconvex function\GPLadd{, and fix the margin $\alpha' \in (0,\alpha)$ of~\cref{subsec:piecewise_approxs}. Since $1-\alpha' > 1-\alpha$, the function $\Psi$ is $(1-\alpha')$-semiconvex as well, and wherever $\Psi$ is twice differentiable there holds $\nabla_{\by}^{2}\Psi(\by) + (1-\alpha')\matr{I}_{n \times n} \succeq (\alpha-\alpha')\matr{I}_{n \times n} \succ 0$; the margin is what keeps the determinant in~\cref{eq:error_mp} below at or above $(\alpha-\alpha')^{n}$}. Following~\cite[Section III]{gaubert2011curse}, we approximate $\Psi$ using $K$ vectors $\{\bp_{k}\}_{k=1}^{K} \subset \Rn$ with $K$ semiconvex functions $\by \mapsto \langle\bp_{k},\by\rangle - \frac{\GPLadd{1-\alpha'}}{2}\normsq{\by}$ and a function $a \colon \Rn \to \R\cup\{+\infty\}$:
\begin{equation}\label{eq:approximation_maxplus}
    \PsiMP(\by) \coloneqq \max_{k \in \{1,\dots,K\}} \left\{\langle\bp_{k},\by\rangle -\frac{\GPLadd{1-\alpha'}}{2}\normsq{\by} - a(\bp_{k})\right\}.
\end{equation}
Here, we suppose the vectors $\{\bp_{k}\}_{k=1}^{K}$ and $\bp \mapsto a(\bp)$ are selected so that $\Psi(\by) \geqslant \PsiMP(\by)$. Let $\mathcal{Y}$ denote a full-dimensional, compact, convex subset of $\Rn$ and consider the $L_{\infty}$ error
\[
\epsilon_{\infty}(\Psi, K, \mathcal{Y}, \PsiMP) \coloneqq \sup_{\by \in \mathcal{Y}} |\Psi(\by) - \PsiMP(\by)|.
\]
Furthermore, we define the corresponding minimal $L_{\infty}$ error as 
\[
\delta_{\infty}(\Psi, K, \mathcal{Y}) = \inf_{\PsiMP \leqslant \Psi} \epsilon_{\infty}(\Psi, K, \mathcal{Y}, \PsiMP).
\]
The following result of~\cite{gaubert2011curse} states that whatever vectors $\{\bp_{k}\}_{k=1}^{K}$ and function $\bp \mapsto a(\bp)$ are used to approximate $\Psi$, the minimal $L_{\infty}$ error scales as an inverse power law in $K$ and the dimension $n$ in the limit $K \to +\infty$.
\begin{theorem}[Gaubert et al. (2011)]\label{thm:maxplus_approx_error}
Let $\alpha > 0$\GPLadd{, let $\alpha' \in (0,\alpha)$}, and let $\mathcal{Y}$ denote a full-dimensional compact, convex subset of $\Rn$. If $\Psi \colon \Rn \to \R$ is twice continuously differentiable and $(1 - \alpha)$\GPLadd{-}semiconvex, then there is a constant $\beta(n) > 0$ such that
\begin{equation}\label{eq:error_mp}
    \delta_{\infty}(\Psi, K, \mathcal{Y})  \sim \beta(n)\left(\frac{1}{K}\int_{\mathcal{Y}} \left(\det\left(\nabla_{\by}^{2}\Psi(\by) + \GPLadd{(1-\alpha')}\matr{I}_{n\times n}\right)\right)^{\frac{1}{2}} \diff\by\right)^{2/n} \ \text{as $K \to +\infty$.}
\end{equation}
\end{theorem}
Thus the minimal $L_{\infty}$ error \GPLadd{scales as $\Theta(1/K^{2/n})$} as $K \to +\infty$, though the rate constant is smaller the closer the Hessian matrix $\nabla_{\by}^{2}\Psi(\by)$ is to the matrix $\GPLadd{-(1-\alpha')\matr{I}_{n \times n}}$\GPLadd{, where the determinant in~\cref{eq:error_mp} degenerates}.

\GPLadd{The exponent $2/n$ in~\cref{eq:error_mp} is \emph{generic} to semiconvex targets. On some priors the max-plus approximation does better or is exact, as the following example shows.}

\begin{example}[\GPLadd{Minimum of quadratics}]\label{ex:quad_exact}
\GPLadd{Let $t > 0$ and let $J = \min_{i \in \{1,\dots,k\}} J_{i}$, where each $J_{i}(\by) = \frac{1}{2}\langle\by,\matr{A}_{i}\by\rangle + \langle\bm{b}_{i},\by\rangle + c_{i}$ is quadratic with $\matr{A}_{i}$ symmetric and $\GPLadd{\matr{I}_{n \times n} + t\matr{A}_{i} \succ 0}$. The Lax--Oleinik semigroup is min-plus linear~\cite{akian2006max,mceneaney2007curse,mceneaney2010convergence}, so there holds $S(\bx,t) = \min_{i} S_{i}(\bx,t)$, where
\[
S_{i}(\bx,t) = \inf_{\by \in \Rn}\left\{\frac{1}{2t}\normsq{\bx-\by} + J_{i}(\by)\right\}.
\]
The objective function above has symmetric positive definite Hessian $\frac{1}{t}(\matr{I}_{n \times n} + t\matr{A}_{i})$, hence strongly convex, and so attains its infimum at the unique stationary point
\[
\by_{i}^{*}(\bx) = \left(\matr{I}_{n\times n} + t\matr{A}_{i}\right)^{-1}\left(\bx - t\bm{b}_{i}\right).
\]
Finally, $\psi(\bx,t) = \frac{1}{2}\normsq{\bx} - t\min_{i} S_{i}(\bx,t) = \max_{i}\left(\frac{1}{2}\normsq{\bx} - tS_{i}(\bx,t)\right)$. The matrices $\matr{A}_{i}$ need not be positive semidefinite: the argument uses only $\matr{I}_{n \times n} + t\matr{A}_{i} \succ 0$, so concave pieces are admitted at small enough $t$.}
\end{example}

\GPLadd{\Cref{ex:quad_exact} suggests that endowing the convex function $\bx \mapsto \psi(\bx,t)$ with structure can help overcome the curse of dimensionality. For instance, let $c > 0$, $k \in \mathbb{N}$, $\{\bm b_{i}\}_{i=1}^{k} \subset \Rn$ and $\{r_{i}\}_{i=1}^{k} \subset \R$. The \emph{max-affine quadratic} (MAQ) function is
\begin{equation}\label{eq:maq_potential}
\psi_{\mathrm{MAQ}}(\by) \coloneqq \frac{c}{2}\normsq{\by} + \max_{i \in \{1,\dots,k\}}\left\{\langle\bm b_{i},\by\rangle + r_{i}\right\}.
\end{equation}
We can interpret the MAQ function as a convex neural network with $1 + k(n+1)$ parameters. Its convex conjugate is given explicitly by
\begin{equation}\label{eq:maq_conjugate}
\psi_{\mathrm{MAQ}}^{*}(\bx) = \min_{\bm\lambda \in \Delta_{k}}\left\{\frac{1}{2c}\normsq{\bx - \textstyle\sum_{i}\lambda_{i}\bm b_{i}} - \sum_{i}\lambda_{i}r_{i}\right\}.
\end{equation}
Thus a MAQ function trained to fit the parameters of~\eqref{eq:maq_potential} through data can recover the learned prior directly via formula~\cref{eq:formula_recovery_prior}.}

\subsection{Applications to the inverse problem for Hamilton--Jacobi Equations}\label{subsec:maxplus_inverse_hjpdes}
We now return to the setting of~\cref{subsec:piecewise_approxs}, in which we learn the backward viscosity solution $\by \mapsto \JBVS(\by)$ from the samples $\{\bx_{k},S(\bx_{k},t),\nabla_{\bx}S(\bx_{k},t)\}_{k=1}^{K}$, assuming additionally that $t\JBVS$ is twice continuously differentiable and $(1-\alpha)$-semiconvex, as in~\cref{cor:jbvs_approx} below. Max-plus algebra theory provides us with the following approximation result.
\begin{corollary}\label{cor:jbvs_approx}
Let $t>0$ and assume $t\JBVS$ is twice continuously differentiable and $(1-\alpha)$-semiconvex with $\alpha > 0$\GPLadd{, with $\alpha'\in(0,\alpha)$ the margin fixed in~\cref{subsec:piecewise_approxs}}. Let $\mathcal{Y}$ denote a full-dimensional compact, convex set of $\Rn$ and let $\tilde{J} \in I_{t}(S)$ denote a function that can reach $\bx \mapsto S(\bx,t)$. Then there exists a constant $\beta(n)$ depending only on $n$ such that
\begin{equation}\label{eq:mp_approx_jbvs}
    \delta_{\infty}(t\JBVS,K,\mathcal{Y}) \sim \beta(n)\left(\frac{1}{K}\int_{\mathcal{Y}} \left(\det\left(t\nabla_{\by}^{2}\JBVS(\by) + \GPLadd{(1-\alpha')}\matr{I}_{n\times n}\right)\right)^{\frac{1}{2}} \diff\by\right)^{2/n}
\end{equation}
as $K \to +\infty$, and
\begin{equation}\label{eq:upperbound_error}
    \delta_{\infty}(\GPLadd{t}\JBVS,K,\mathcal{Y}) \leqslant t\sup_{\by \in \mathcal{Y}}|\tilde{J}(\by) - \JPQM(\by)|.
\end{equation}
\end{corollary}
\begin{proof}
\GPLadd{The asymptotic~\cref{eq:mp_approx_jbvs} is~\cref{thm:maxplus_approx_error} applied to $\Psi \equiv t\JBVS$, whose hypotheses hold by assumption. For~\cref{eq:upperbound_error}, note that $t\JPQM$ is admissible in the infimum defining $\delta_{\infty}(t\JBVS,K,\mathcal{Y})$. Indeed, by~\cref{eq:yk_formula} and Fenchel duality, we have $\bx_{k} \in \partial\left(t\JBVS + \frac{1}{2}\normsq{\cdot}\right)(\by_{k})$, so expanding the square in the definition of $\JPQM$ puts $t\JPQM$ in the form~\cref{eq:approximation_maxplus} with $\bp_{k} = \bx_{k} - \alpha'\by_{k}$, and $\JPQM \leqslant \JBVS$ by~\cref{subsec:piecewise_approxs}. Since~\cref{thm:prior_optimality} gives $\JPQM \leqslant \JBVS \leqslant \tilde{J}$ on $\Rn$, there holds
\[
\delta_{\infty}(t\JBVS,K,\mathcal{Y}) \leqslant t\sup_{\by \in \mathcal{Y}}\left(\JBVS(\by) - \JPQM(\by)\right) \leqslant t\sup_{\by \in \mathcal{Y}}\left|\tilde{J}(\by) - \JPQM(\by)\right|.
\]}
\end{proof}
\GPLadd{Thus the best minorant of the form~\cref{eq:approximation_maxplus} approaches $t\JBVS$ at the rate $\Theta(1/K^{2/n})$, and no reachable $\tilde{J} \in I_{t}(S)$ lies closer to $\JPQM$ than $\delta_{\infty}(t\JBVS,K,\mathcal{Y})/t$ on $\mathcal{Y}$. The rate constant is small when the eigenvalues of $t\nabla_{\by}^{2}\JBVS$ cluster near $-(1-\alpha')$, where the determinant in~\cref{eq:mp_approx_jbvs} degenerates.}

\GPLadd{\Cref{eq:upperbound_error} does not bound the error of $\JPQM$. The sample spacing controls the latter: if $\nabla_{\by}^{2}\left(t\JBVS + \frac{1}{2}\normsq{\cdot}\right) \preceq L\matr{I}_{n\times n}$ on $\mathcal{Y}$ for some $L>0$, then with $h_{K}(\mathcal{Y}) \coloneqq \sup_{\by \in \mathcal{Y}}\min_{k}\norm{\by-\by_{k}}$, the descent lemma applied at each $\by_{k}$ gives
\begin{equation}\label{eq:pqm_fill}
\sup_{\by \in \mathcal{Y}}\left(\JBVS(\by) - \JPQM(\by)\right) \leqslant \frac{L-\alpha'}{2t}\,h_{K}(\mathcal{Y})^{2},
\end{equation}
while covering $\mathcal{Y}$ by the balls $B(\by_{k},h_{K}(\mathcal{Y}))$ forces $h_{K}(\mathcal{Y}) \geqslant \left(\operatorname{vol}\mathcal{Y}/\omega_{n}K\right)^{1/n}$, with $\omega_{n}$ the volume of the unit ball. For well-spread samples $h_{K}(\mathcal{Y}) = \Theta(K^{-1/n})$, so the right-hand side of~\cref{eq:pqm_fill} is $O(K^{-2/n})$. The matching lower bound is~\cref{eq:mp_approx_jbvs} itself: $t\JPQM$ is admissible in the infimum defining $\delta_{\infty}(t\JBVS,K,\mathcal{Y})$, so the left-hand side of~\cref{eq:pqm_fill} is at least $\delta_{\infty}(t\JBVS,K,\mathcal{Y})/t$, which is $\Theta(K^{-2/n})$. The minorant built from the $K$ samples therefore errs by $\Theta(K^{-2/n})$ and attains the optimal max-plus exponent.}

\GPLadd{In practice we fit a convex neural network $\hat\psi$ to $\psi(\cdot,t)$ rather than use the minorant $\JPQM$ that~\cref{cor:jbvs_approx} bounds. The next result converts an error on $\hat\psi$ into an error on the recovered prior at the price $1/t$, provided the two points at which the conjugates are attained lie in the region where $\hat\psi$ was fitted.}

\begin{proposition}[\GPLadd{Transfer of the potential error to the prior}]\label{prop:transfer_values}
\GPLadd{Let $t>0$, assume $\bx \mapsto S(\bx,t)$ is reachable, and let $B \subset \Rn$ and $\hat\psi \colon \Rn \to \R$ satisfy $\varepsilon_{\psi}(t) \coloneqq \sup_{\bx \in B}\left|\hat\psi(\bx) - \psi(\bx,t)\right| < +\infty$. Fix $\by \in \Rn$ and suppose that $\bx \mapsto \langle\bx,\by\rangle - \psi(\bx,t)$ attains its supremum over $\Rn$ at some $\bx^{*} \in B$, and that $\bx \mapsto \langle\bx,\by\rangle - \hat\psi(\bx)$ attains its supremum at some $\hat\bx \in B$. Then the prior recovered from $\hat\psi$ by conjugation, $\hat{J} \coloneqq \frac{1}{t}\left(\hat\psi^{*} - \frac{1}{2}\normsq{\cdot}\right)$, satisfies
\begin{equation}\label{eq:pointwise_bound}
\left|\hat{J}(\by) - \JBVS(\by)\right| \leqslant \frac{1}{t}\max\left\{\left|\hat\psi-\psi(\cdot,t)\right|(\bx^{*}),\ \left|\hat\psi-\psi(\cdot,t)\right|(\hat\bx)\right\} \leqslant \frac{\varepsilon_{\psi}(t)}{t}.
\end{equation}}
\end{proposition}

\begin{proof}
\GPLadd{Since $\hat\bx \in B$, there holds
\[
\hat\psi^{*}(\by) = \langle\hat\bx,\by\rangle - \hat\psi(\hat\bx) \leqslant \langle\hat\bx,\by\rangle - \psi(\hat\bx,t) + \left|\hat\psi-\psi(\cdot,t)\right|(\hat\bx) \leqslant \psi^{*}(\by,t) + \left|\hat\psi-\psi(\cdot,t)\right|(\hat\bx),
\]
and exchanging the roles of $\psi(\cdot,t)$ and $\hat\psi$ at $\bx^{*}$ gives the reverse inequality. Subtracting $\frac{1}{2}\normsq{\by}$ from both conjugates and dividing by $t$ yields~\cref{eq:pointwise_bound} through~\cref{eq:formula_recovery_prior}.}
\end{proof}

\begin{remark}[\GPLadd{Why the training region must cover the preimages}]\label{rem:box}
\GPLadd{The point $\bx^{*}$ is a preimage of $\by$ under the proximal map $\nabla_{\bx}\psi(\cdot,t)$, and $\hat\bx$ is the point returned when $\hat\psi^{*}(\by)$ is evaluated. A margin is sufficient a priori because if $\hat\psi$ is convex and $\psi(\cdot,t)$ is $\sigma$-strongly convex near $B$, then $\operatorname{dist}(\bx^{*},\partial B) > 2\sqrt{\varepsilon_{\psi}/\sigma}$ places $\hat\bx$ in $B$ as well (since $\psi(\cdot,t) - \langle\cdot,\by\rangle$ then exceeds its minimum by more than $2\varepsilon_{\psi}$ on $\partial B$).}
\end{remark}


\section{Numerical results}\label{sec:numerics}
In this section, we leverage the theory of viscosity solutions of HJ PDEs to learn proximal operators of certain (non)convex priors. \GPLadd{To do this, we use the learned proximal network (LPN) architecture of~\cite{fang2024whats} and the max-affine quadratic function~\eqref{eq:maq_potential} discussed in~\cref{subsec:maxplus}. We present our setup in~\cref{subsec:setup}. In \cref{subsec:closed_form}, we report numerical experiments on four priors with closed-form proximal mappings, namely the $\ell_{1}$ norm, a minimum of two quadratics, a concave quadratic, and the negative $\ell_{1}$ norm. Following the work of~\cite{darbon2021bayesian}, who showed that a total variation posterior mean denoiser is the proximal mapping of some implicit prior, we show numerically in~\cref{subsec:tv_exp} that the implicit prior can be successfully learned from data. Finally, in~\cref{subsec:pnp_exp} we recover priors from denoisers on two problems with exact references, and then diagnose and recover the regularizer of the two plug-and-play operators released in~\cite{hu2023plug}.}

\subsection{\GPLadd{Experimental setup}}\label{subsec:setup}

\subsubsection{\GPLadd{Networks, data and training}}\label{subsubsec:protocol}
\GPLadd{We fix $t = 1$ in~\cref{subsec:closed_form} and give $t$ with each later experiment. We write $d$ for the dimension $n$ when it indexes an experiment.}

\GPLadd{\emph{Networks.} We use two convex neural networks. The first is the learned proximal network (LPN) of~\cite{fang2024whats}, a fully connected input-convex neural network with Softplus activation $\mathrm{softplus}_{\beta}(s) = \beta^{-1}\ln(1+e^{\beta s})$ at $\beta = 5$, which we take with $2$ hidden layers of $256$ units. After every optimizer step, we clip to zero the negative entries of the weight matrices that carry the hidden features forward, which by~\cite[Prop.~3.1]{fang2024whats} keeps the network convex in its input. The width is $256$ at every $d$, so that $d$ is the only variable across a row. The parameter count $P$ then varies only through the input layer, from $133\,379$ at $d = 2$ to $181\,057$ at $d = 64$, while the sample size $N$ grows proportionally to $d$. The ratio $N/P$ runs from $0.22$ at $d = 2$ to $5.3$ at $d = 64$, so the low-dimensional runs are the over-parameterized ones.}

\GPLadd{The second class is the max-affine quadratic (MAQ) function~\cref{eq:maq_potential} of~\cref{subsec:maxplus}, which we read as a convex network with parameters $c > 0$, $\{\bm b_{i}\}_{i=1}^{k}$ and $\{r_{i}\}_{i=1}^{k}$. We write $c$ as a softplus of an unconstrained scalar to keep it positive, and convexity then needs no clipping, since a maximum of affine functions plus a positive multiple of $\normsq{\by}$ is convex by construction. We use $k = 2$ atoms at every dimension and on every family, which gives $1 + k(d+1) = 2d+3$ parameters, from $7$ at $d = 2$ to $131$ at $d = 64$, three to four orders of magnitude below the LPN count. The conjugate~\cref{eq:maq_conjugate} is available in closed form, so a fitted MAQ network returns the prior with no second network and no inversion.}

\GPLadd{\emph{Data.} We report all errors on the \emph{query box} $[-4,4]^{d}$, and train on a possibly larger \emph{training box} $[-A,A]^{d}$, with $A$ the smallest half-width for which the image of $[-A,A]^{d}$ under $\nabla_{\bx}\psi(\cdot,t)$ contains the query box. That gives $A = 5$ for $\normone{\bx}$, $A = 9$ for the min-plus prior, and $A = 4$ for the two priors whose proximal map contracts. Every recovery below evaluates $\psi$ at preimages of the query points, so a training box failing to cover them would violate the hypothesis of~\cref{prop:transfer_values} at the very points where we measure the error. The training set is $N = 15\,000\,d$ points drawn uniformly from the training box with seed $1$, the validation set $4000$ fresh points with seed $2$, and the test set $4000$ points on the query box with seed $3$, of which the first $1000$ are shared by every method in every table.}

\GPLadd{\emph{Training.} We train each LPN with the Adam optimizer, mini-batches of $512$, a budget of $S = 2.5\times10^{5}$ optimizer \emph{steps}, and learning rate $10^{-3}$ multiplied by $10^{-1}$ at $50\%$ and at $75\%$ of $S$. We use no early stopping and keep the checkpoint with the lowest validation error, so the validation error we report describes the LPN we return. Every LPN minimizes a mean-squared error (MSE). We fit the MAQ network differently, by full-batch Adam on $8000$ of the training points, $4000$ steps at learning rate $5\times10^{-2}$ multiplied by $10^{-1}$ at half the budget, again minimizing the MSE against samples of $\psi(\cdot,t)$. All runs are on CPU.}

\GPLadd{\emph{Recovering the priors.} We compare four ways of recovering priors from data.}
\begin{enumerate}
\item[(i)] \GPLadd{The \emph{iterative} recovery is the regularized inversion of~\cite{fang2024whats}. For each query point $\bx$ it minimizes $\by \mapsto \psiNN(\by) + \frac{\alpha_{\mathrm{inv}}}{2}\normsq{\by} - \langle\bx,\by\rangle$, one optimization per point, and sets $J(\bx) = \langle\bx,\by\rangle - \psiNN(\by) - \frac{1}{2}\normsq{\bx}$ at the minimizer. Its only knob, $\alpha_{\mathrm{inv}} \geqslant 0$, is not free, since at the solution $\nabla\psiNN(\by) = \bx - \alpha_{\mathrm{inv}}\by$, so the prior is biased at order $\alpha_{\mathrm{inv}}$, while at $\alpha_{\mathrm{inv}} = 0$ the objective need not be coercive and the iterate may diverge. We therefore run every case of~\cref{subsec:closed_form} at $\alpha_{\mathrm{inv}} \in \{0, 0.1\}$ and report the better of the two, recording the winner in the tables as $\alpha_{\mathrm{inv}}^{*}$. The experiments of~\cref{subsec:tv_exp,subsec:pnp_exp} run it unregularized, at $\alpha_{\mathrm{inv}} = 0$.}
\item[(ii)] \GPLadd{The \emph{two-network} recovery is the two-phase procedure of~\cref{subsubsec:two_phase}, in which a second network $\JNN$ trains on the pairs $\{\by_{k},G_{k}\}$ that $\psiNN$ produces by conjugation and then evaluates the prior in one forward pass.}
\item[(iii)] \GPLadd{The \emph{one-network} recovery trains a single network on data the proximal map supplies. Where the value function is available at the samples, we form the pairs $\{\by_{k},G_{k}\}$ directly, as in~\cref{subsec:closed_form}. Where only the map is available we train $J_{\theta}$ against the optimality condition $\nabla J_{\theta}(\by_{k}) = \bx_{k} - \by_{k}$, as in~\cref{subsec:tv_exp,subsec:pnp_exp}.}
\item[(iv)] \GPLadd{The \emph{MAQ} recovery trains no LPN for the prior. It fits the MAQ network to samples of $\psi(\cdot,t)$ and recovers the prior by the conjugation~\cref{eq:maq_conjugate}.}
\end{enumerate}

\GPLadd{Recoveries (i), (ii), (iv) and the value-supervised form of (iii) return a convex function minus $\frac{1}{2}\normsq{\cdot}$ and are therefore $1$-semiconvex. The gradient-supervised form of (iii) returns $J_{\theta}$ itself, which the LPN architecture makes convex, so that form assumes that the prior is convex, whereas the others do not.}

\GPLadd{\emph{Error measures.} For a recovered prior $\hat{J}$ and a reference $J_{\mathrm{ref}}$ evaluated at the $1000$ shared test points $\{\bz_{k}\}$, the \emph{relative $L^{2}$ error} is
\begin{equation}\label{eq:relL2}
\mathrm{relL}^{2}(\hat{J},J_{\mathrm{ref}}) \coloneqq \left(\sum_{k=1}^{1000}\bigl|\hat{J}(\bz_{k}) - J_{\mathrm{ref}}(\bz_{k})\bigr|^{2}\right)^{1/2}\Big/\left(\sum_{k=1}^{1000}\bigl|J_{\mathrm{ref}}(\bz_{k})\bigr|^{2}\right)^{1/2},
\end{equation}
and the root-mean-squared error (RMSE) is the numerator divided by $\sqrt{1000}$. We take $J_{\mathrm{ref}} = \JBVS$, because~\cref{thm:prior_optimality} identifies the backward viscosity solution $\JBVS$ as the prior to recover.\footnote{\GPLadd{The min-plus and concave priors are not uniformly Lipschitz continuous, so they lie outside the class $I_{t}(S)$ of~\cref{eq:reachable_init_conds} on which~\cref{thm:prior_optimality} is stated. The identification survives on both without viscosity theory, since~\cref{eq:formula_recovery_prior} is convex duality and asks only that $\bx \mapsto \psi(\bx,t)$ be convex, which we check in~\cref{subsubsec:nonconvex,subsubsec:concave}.}} The two coincide for the $\ell_{1}$ and concave priors below, and differ for $-\normone{\bx}$ by a few percent and for the min-plus prior by a few parts in $10^{4}$. On both we report the separation alongside the recovery, so the reader can see when a method resolves the two references. \emph{Validation MSE} means the mean-squared error of a network against its own training target on the $4000$ validation points of seed $2$, never on training data. It measures the fit of one network, while the relative $L^{2}$ error measures the recovered prior and is what we compare across methods.}

\subsubsection{\GPLadd{The two-phase procedure}}\label{subsubsec:two_phase}
\GPLadd{At every sample we evaluate the optimal value $S(\bx_{k},1)$ of the proximal problem. In the first phase we learn the value function $S(\bx,1)$ of the HJ equation~\cref{eqn:intro2} with initial data $J$ at $t = 1$, encoded in a convex potential $\psiNN(\by;\theta)$ trained by minimizing the MSE loss}
\[
\min_{\theta}\frac{1}{N} \sum_{k=1}^{N} \left(\psiNN(\bx_{k};\theta) - \left(\frac{1}{2}\normsq{\bx_{k}}- S(\bx_{k},1)\right)\right)^2.
\]
\GPLadd{We then read the value function off as} $\bx \mapsto \frac{1}{2}\normsq{\bx} - \psiNN(\bx;\theta)$.

\GPLadd{In the second phase we train a separate LPN} $\by \mapsto \JNN(\by;\tilde{\theta})$ \GPLadd{on the conjugate $\psi^{*}(\cdot,1)$, from which~\cref{eq:formula_recovery_prior} returns the prior $\JBVS$ of the terminal value HJ PDE~\cref{eq:finite_time_backward}. Formulas~\cref{eq:minimum_representation} and~\cref{eq:formula_recovery_prior} give} approximate ``samples'' $\{\by_{k},G_{k}\}_{k=1}^{N}$ \GPLadd{of that conjugate}:
\[
\by_{k} = \nabla_{\bx}\psiNN(\bx_{k};\theta), \quad \GPLadd{\psi^{*}(\by_{k},1)} \approx G_{k} \coloneqq \langle\GPLadd{\bx_{k}},\by_{k}\rangle - \psiNN(\GPLadd{\bx_{k};\theta}),
\]
\GPLadd{and we train} the second LPN by minimizing the MSE loss
\[
\min_{\tilde{\theta}} \frac{1}{N} \sum_{k=1}^{N} \left(\JNN(\by_{k};\tilde{\theta}) - G_{k}\right)^{2}.
\]
\GPLadd{We read the prior off as} $\by \mapsto \JNN(\by;\tilde{\theta}) - \frac{1}{2}\normsq{\by}$ \GPLadd{and compare it against the reference $\JBVS$ of the experiment at hand. Fitting the conjugate keeps the network's target convex, and the subtraction is what makes the recovered prior $1$-semiconvex.}

\subsection{Priors known in closed form}\label{subsec:closed_form}
\subsubsection{Convex prior: the $\ell_{1}$ norm}\label{subsubsec:l1}
\begin{figure}[!htbp]
 \centering
\includegraphics[width=0.8\textwidth,height=0.40\textheight]{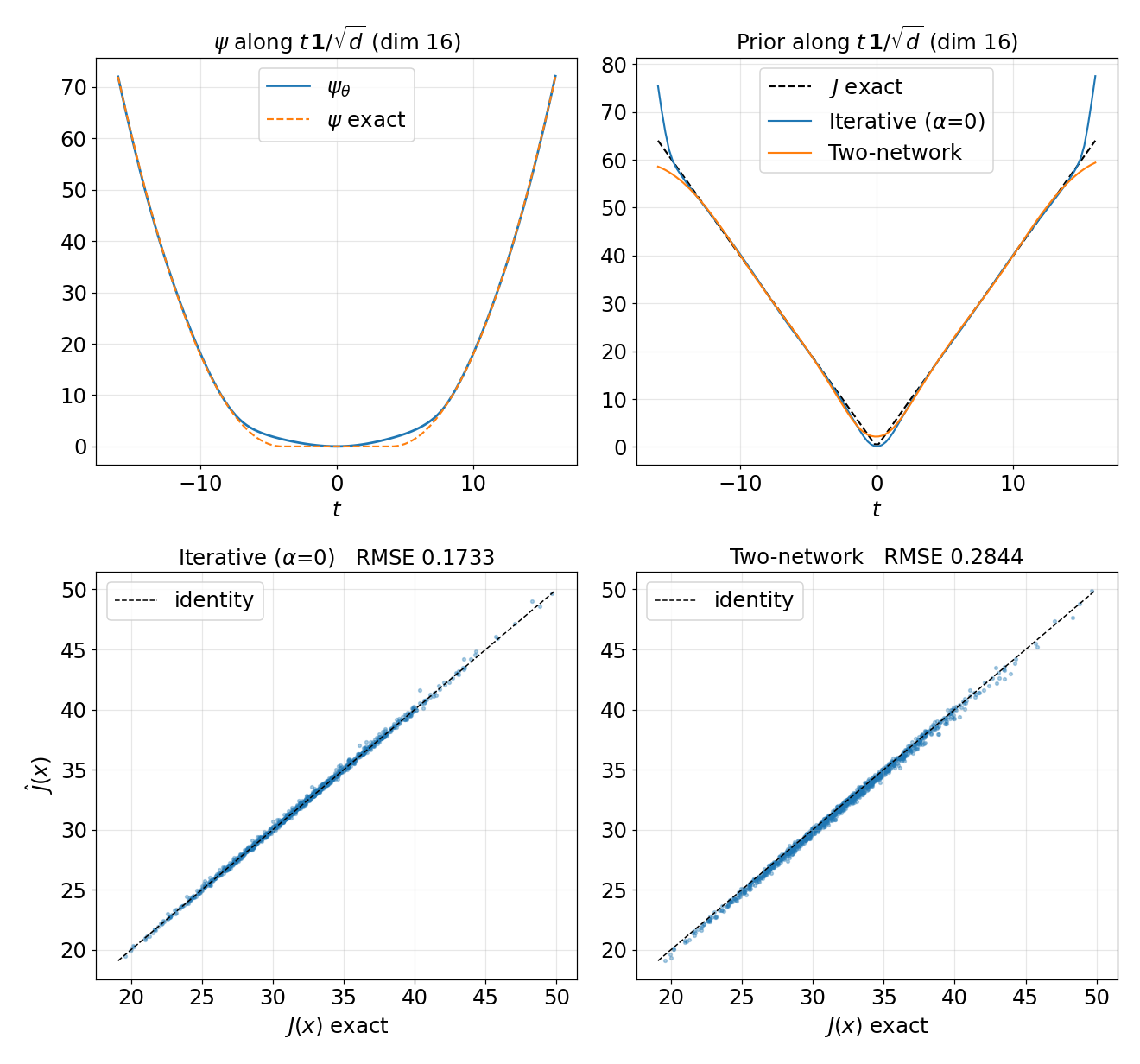}
\caption{\GPLadd{Convex $\ell_1$ prior, dimension $16$. Left: cross sections of the convex potential $\psiNN(\bx;\theta)$ and of the prior recovered by both methods, against the exact $\bx \mapsto \frac{1}{2}\normsq{\bx} - S(\bx,1)$ and against $\JBVS$. Right: the root-mean-squared error of the recovered prior against $\JBVS$, over the shared test points. The iterative recovery of~\cite{fang2024whats} inverts $\psiNN$ once per query point; the two-network recovery $\JNN(\by;\tilde{\theta})$ uses a single forward pass.}}
\label{dim18}
\end{figure}

We begin with the convex prior $J(\bx) = \normone{\bx}$. \GPLadd{The} value function $S(\bx,t)$ \GPLadd{follows from} the soft-thresholding operator~\cref{eq:soft_thresholding}, and $\JBVS = J$. \GPLadd{\Cref{F1} reports the recovery,~\cref{F1b} the network fits, and~\cref{dim18} plots both at $d = 16$.}

\GPLadd{Both networks fit $\psi(\cdot,1)$ and the prior well through $d = 16$ and lose accuracy at $d = 32$ and $64$. As Table 2 shows, the validation MSE of $\psiNN$ rises from $1.17\times10^{-2}$ at $d = 16$ to $15.2$ at $d = 64$. The recovered prior degrades less than that fit suggests, because the relative $L^{2}$ error normalizes by a reference whose magnitude also grows with $d$. The one-network variant is the most stable in the dimension, with its error growing by a factor of $5.8$ against a factor of $65$ (two-network) and a factor of 118 (iterative inversion) from $d = 2$ to $d = 64$. The MAQ recovery is weakest here because it does not apply to $\normone{\bx}$, whose semiconvex transform is a maximum of $2^{d}$ affine pieces. Thus~\cref{F1} shows that the one-network approach works significantly better in high dimensions, while the iterative inversion works well only in lower dimensions.}

\begin{table}[!htbp]
\centering
\caption{\GPLadd{Convex $\ell_{1}$ prior: relative $L^{2}$ error~\cref{eq:relL2} of the recovered prior against $\JBVS = J$, on the $1000$ shared test points. Here $\alpha_{\mathrm{inv}}^{*}$ is the regularization parameter chosen by the best-of-$\{0,0.1\}$ rule. Boldface marks the best entry in a row.}}
\footnotesize
\setlength{\tabcolsep}{5pt}
\begin{tabular}{c cccc}
\toprule
$d$ & iterative ($\alpha_{\mathrm{inv}}^{*}$) & two-network & one-network & MAQ\\
\midrule
2 & $\bm{1.53\times 10^{-3}}$ (0) & $3.18\times 10^{-3}$ & $3.77\times 10^{-3}$ & $1.25\times 10^{-1}$\\
4 & $\bm{1.54\times 10^{-3}}$ (0) & $3.78\times 10^{-3}$ & $4.24\times 10^{-3}$ & $9.52\times 10^{-2}$\\
8 & $\bm{1.69\times 10^{-3}}$ (0) & $4.66\times 10^{-3}$ & $5.40\times 10^{-3}$ & $7.50\times 10^{-2}$\\
16 & $1.10\times 10^{-2}$ (0) & $1.07\times 10^{-2}$ & $\bm{9.77\times 10^{-3}}$ & $5.66\times 10^{-2}$\\
32 & $4.90\times 10^{-2}$ (0.1) & $7.68\times 10^{-2}$ & $\bm{1.48\times 10^{-2}}$ & $5.21\times 10^{-2}$\\
64 & $1.81\times 10^{-1}$ (0.1) & $2.06\times 10^{-1}$ & $\bm{2.20\times 10^{-2}}$ & $1.28\times 10^{-1}$\\
\bottomrule
\end{tabular}
\label{F1}
\end{table}

\begin{table}[!htbp]
\centering
\caption{\GPLadd{Convex $\ell_{1}$ prior: validation mean-squared error of each trained network against its own target, on the $4000$ validation points of seed $2$.}}
\footnotesize
\setlength{\tabcolsep}{4pt}
\begin{tabular}{l cccccc}
\toprule
network & $d = 2$ & $4$ & $8$ & $16$ & $32$ & $64$\\
\midrule
$\psiNN$ & $5.81\times 10^{-5}$ & $2.21\times 10^{-4}$ & $8.27\times 10^{-4}$ & $1.17\times 10^{-2}$ & $4.10\times 10^{0}$ & $1.52\times 10^{1}$\\
$\JNN$ (two-net) & $1.60\times 10^{-4}$ & $6.42\times 10^{-4}$ & $3.39\times 10^{-3}$ & $4.45\times 10^{-2}$ & $1.64\times 10^{0}$ & $1.31\times 10^{1}$\\
$\JNN$ (one-net) & $2.48\times 10^{-4}$ & $1.17\times 10^{-3}$ & $5.45\times 10^{-3}$ & $5.10\times 10^{-2}$ & $4.13\times 10^{-1}$ & $5.21\times 10^{0}$\\
\bottomrule
\end{tabular}
\label{F1b}
\end{table}

\subsubsection{\GPLadd{A non-convex prior from min-plus algebra}}\label{subsubsec:nonconvex}
\begin{figure}[!htbp]
 \centering
\includegraphics[width=0.8\textwidth,height=0.40\textheight]{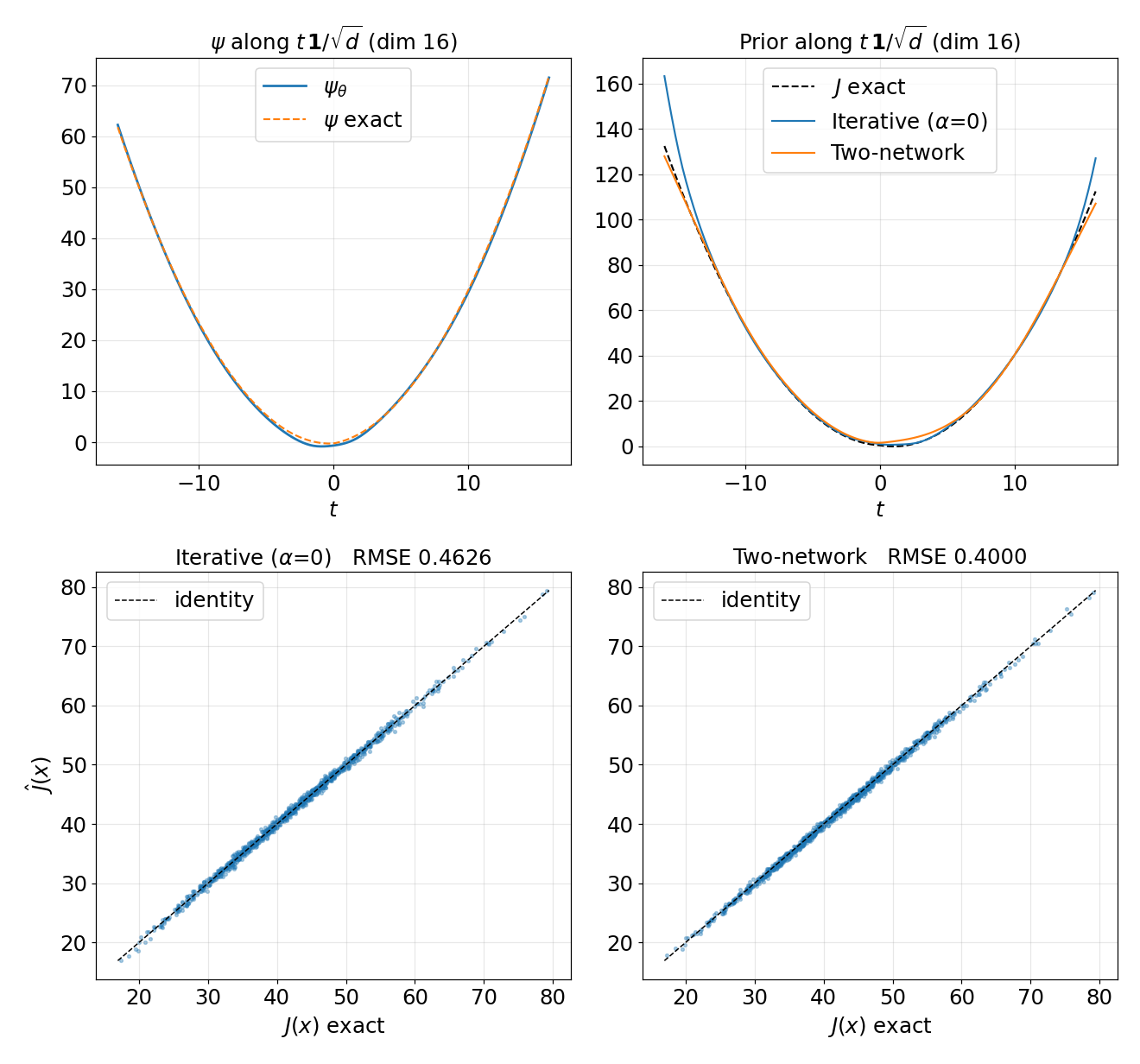}
\caption{\GPLadd{Min-plus prior, dimension $16$. Left: cross sections of the convex potential $\psiNN(\bx;\theta)$ against the exact $\bx \mapsto \frac{1}{2}\normsq{\bx}-S(\bx,1)$, and of the recovered prior against $\JBVS$. Right: the root-mean-squared error of the recovered prior against $\JBVS$ on the shared test points.}}
\label{dim28}
\end{figure}

\begin{figure}[!htbp]
 \centering
\includegraphics[width=0.8\textwidth,height=0.40\textheight]{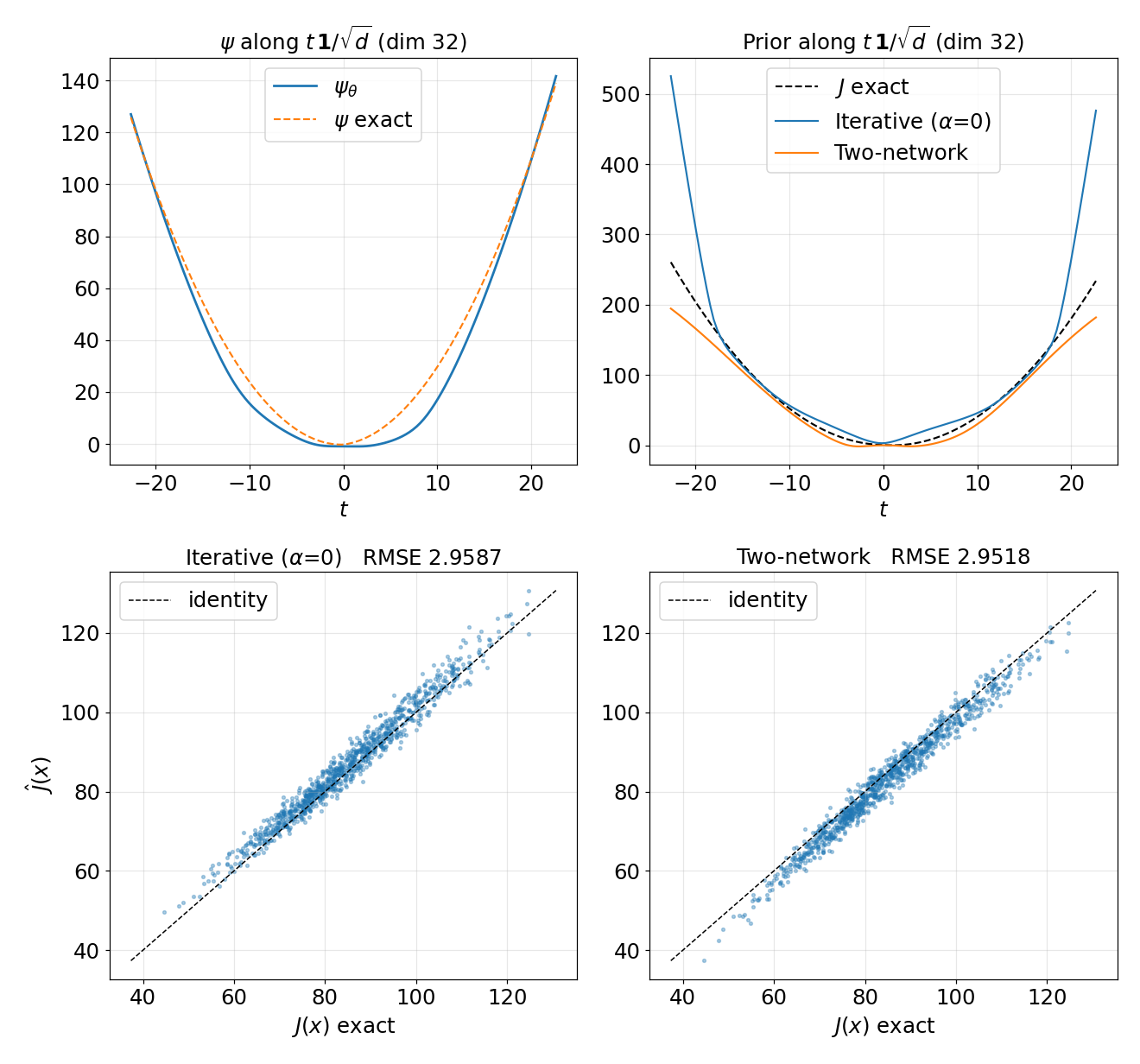}
\caption{\GPLadd{Min-plus prior, dimension $32$. Panels as in~\cref{dim28}.}}
\label{dim282}
\end{figure}

\begin{figure}[h!]
\centering
\includegraphics[width=0.8\textwidth,height=0.40\textheight]{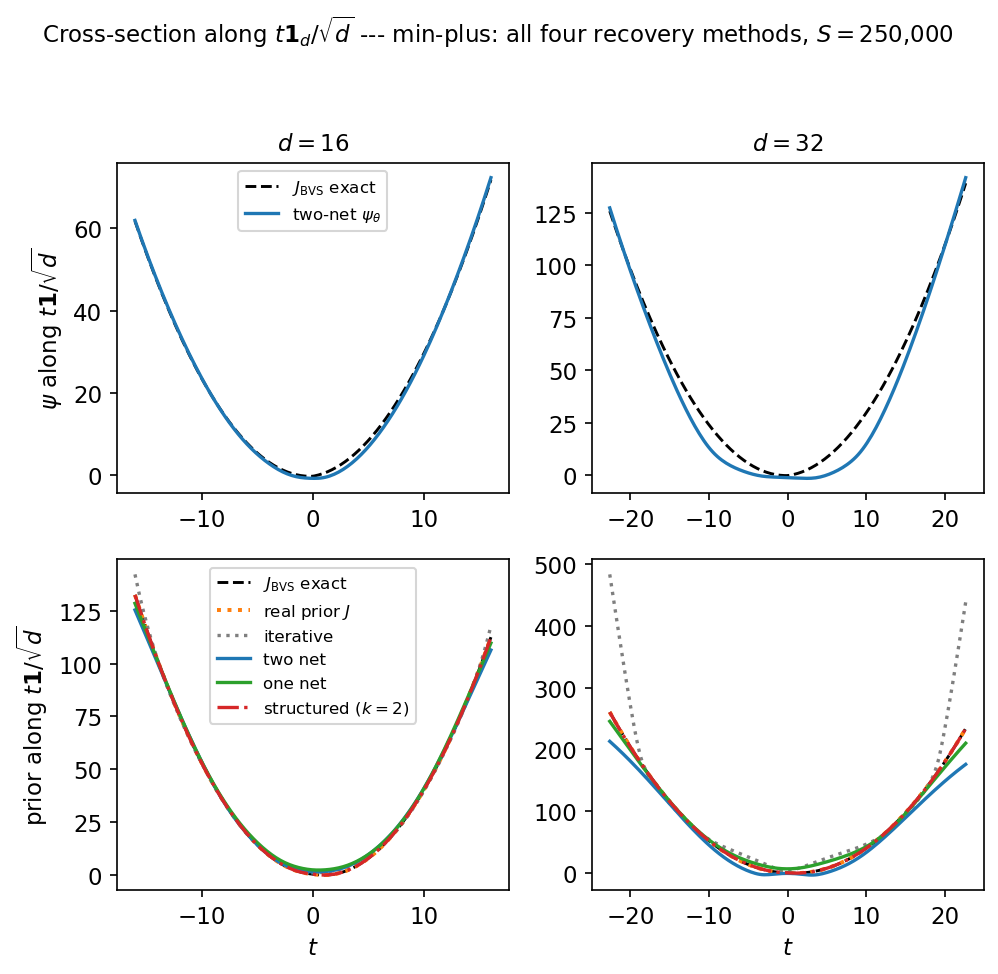}
\caption{\GPLadd{Min-plus prior, cross sections along the diagonal ray $\bx(s) = s\boldsymbol{1}_{d}/\sqrt{d}$ at $d = 16$ and $d = 32$. Top: the convex potential $\psi(\cdot,1)$, exact and as recovered. Bottom: the recovered prior for all four recoveries, plotted against both references $\JBVS$ and $J$. The MAQ recovery is indistinguishable from $\JBVS$. The iterative recovery diverges near the ends of the ray, where the preimages leave the region on which $\psiNN$ was trained.}}
\label{fig:ray_minplus}
\end{figure}

We consider the min-plus prior
\[
J(\bx) = \min\left(\frac{1}{2\sigma_{1}}\normsq{\bx-\mu_{1}}, \frac{1}{2\sigma_{2}}\normsq{\bx - \mu_{2}}\right),
\]
with $\mu_{1} = (1,0,\dots,0)$, $\mu_{2} = \boldsymbol{1}/\sqrt{n}$, and $\sigma_{1} = \sigma_{2} = 1$. \GPLadd{This prior is not uniformly Lipschitz continuous, the assumption under which the HJ PDE has a unique viscosity solution. We include it to see whether an LPN learns it regardless. The identification $J_{\mathrm{ref}} = \JBVS$ still holds, by the duality of~\cref{eq:formula_recovery_prior}: at $\sigma_{1} = \sigma_{2} = 1$ the potential $\psi(\bx,1) = \max_{i}\left(\frac{1}{2}\normsq{\bx} - \frac{1}{4}\normsq{\bx-\mu_{i}}\right)$ has pieces of Hessian $\frac{1}{2}\matr{I}_{d \times d} \succ 0$, hence is convex.} \GPLadd{\Cref{ex:quad_exact} applies here because $J$ is a minimum of two quadratics with $\matr{I}_{n \times n} + t\matr{A}_{i} = 2\matr{I}_{n \times n} \succ 0$, hence $\psi(\cdot,1)$ is an exact maximum of two quadratics and the MAQ potential~\cref{eq:maq_potential} with $k = 2$ represents it exactly at any $d$. In~\cref{F2}, MAQ is accurate to $10^{-8}$ at $d = 2$ and to $10^{-4}$ at $d = 64$, beating the best of the three networks by a factor of $249$ at $d = 16$ and $365$ at $d = 64$, and its error stays nearly flat in the dimension where every network error grows. \Cref{F2b} reports the network fits. \Cref{dim28,dim282} plot cross sections at $d = 16$ and $d = 32$, and~\cref{fig:ray_minplus} follows a diagonal ray, near whose ends the iterative recovery diverges because its preimages leave the region on which $\psiNN$ was trained. As Table~\cref{F2} shows, the MAQ neural network significantly outperforms all other trained neural networks and scale well to high dimensions.} 

\begin{table}[!htbp]
\centering
\caption{\GPLadd{Min-plus prior: relative $L^{2}$ error of the recovered prior against $\JBVS$, columns as in~\cref{F1}. The MAQ potential with $k = 2$ atoms represents $\psi(\cdot,1)$ exactly, so its error reflects only the fit of the two atoms to sampled data.}}
\footnotesize
\setlength{\tabcolsep}{5pt}
\begin{tabular}{c cccc}
\toprule
$d$ & iterative ($\alpha_{\mathrm{inv}}^{*}$) & two-network & one-network & MAQ\\
\midrule
2 & $1.48\times 10^{-3}$ (0) & $1.20\times 10^{-3}$ & $1.32\times 10^{-3}$ & $\bm{1.52\times 10^{-8}}$\\
4 & $2.37\times 10^{-3}$ (0) & $2.04\times 10^{-3}$ & $2.11\times 10^{-3}$ & $\bm{6.91\times 10^{-9}}$\\
8 & $4.11\times 10^{-3}$ (0) & $3.71\times 10^{-3}$ & $3.46\times 10^{-3}$ & $\bm{1.10\times 10^{-7}}$\\
16 & $1.09\times 10^{-2}$ (0) & $1.00\times 10^{-2}$ & $8.37\times 10^{-3}$ & $\bm{3.36\times 10^{-5}}$\\
32 & $3.57\times 10^{-2}$ (0) & $3.64\times 10^{-2}$ & $2.05\times 10^{-2}$ & $\bm{1.53\times 10^{-4}}$\\
64 & $5.35\times 10^{-2}$ (0.1) & $7.89\times 10^{-2}$ & $4.43\times 10^{-2}$ & $\bm{1.21\times 10^{-4}}$\\
\bottomrule
\end{tabular}
\label{F2}
\end{table}

\begin{table}[!htbp]
\centering
\caption{\GPLadd{Min-plus prior: validation mean-squared error, as in~\cref{F1b}.}}
\footnotesize
\setlength{\tabcolsep}{4pt}
\begin{tabular}{l cccccc}
\toprule
network & $d = 2$ & $4$ & $8$ & $16$ & $32$ & $64$\\
\midrule
$\psiNN$ & $2.16\times 10^{-5}$ & $5.31\times 10^{-4}$ & $9.57\times 10^{-3}$ & $3.13\times 10^{-1}$ & $8.72\times 10^{0}$ & $6.55\times 10^{1}$\\
$\JNN$ (two-net) & $1.27\times 10^{-5}$ & $5.87\times 10^{-4}$ & $1.54\times 10^{-2}$ & $5.14\times 10^{-1}$ & $1.19\times 10^{1}$ & $5.86\times 10^{1}$\\
$\JNN$ (one-net) & $1.50\times 10^{-5}$ & $2.50\times 10^{-4}$ & $6.35\times 10^{-3}$ & $1.91\times 10^{-1}$ & $3.85\times 10^{0}$ & $6.73\times 10^{1}$\\
\bottomrule
\end{tabular}
\label{F2b}
\end{table}

\subsubsection{Concave prior}\label{subsubsec:concave}
We consider the concave prior $J(\bx) = -\normsq{\bx}/4$. \GPLadd{It is not uniformly Lipschitz continuous and therefore lies outside the class $I_{t}(S)$ on which~\cref{thm:prior_optimality} is stated. The identification of $\JBVS$ below rests on convex duality instead.} \GPLadd{Here} $S(\bx,t) = -\normsq{\bx}/(2(2-t)$ for $0 \leqslant t < 2$, so $S(\bx,1) = -\frac{1}{2}\normsq{\bx}$ and $\JBVS = J$\GPLadd{: the potential $\psi(\bx,1) = \frac{1}{2}\normsq{\bx} - S(\bx,1) = \normsq{\bx}$ is convex with conjugate $\psi^{*}(\by,1) = \frac{1}{4}\normsq{\by}$, and~\cref{eq:formula_recovery_prior} then returns $\JBVS(\by) = \frac{1}{4}\normsq{\by} - \frac{1}{2}\normsq{\by} = J(\by)$}. \GPLadd{The prior is a single quadratic with $\matr{I}_{n\times n} + t\matr{A} = \frac{1}{2}\matr{I}_{n\times n} \succ 0$, so~\cref{ex:quad_exact} applies here.}

\begin{figure}[!htbp]
 \centering
\includegraphics[width=0.8\textwidth,height=0.40\textheight]{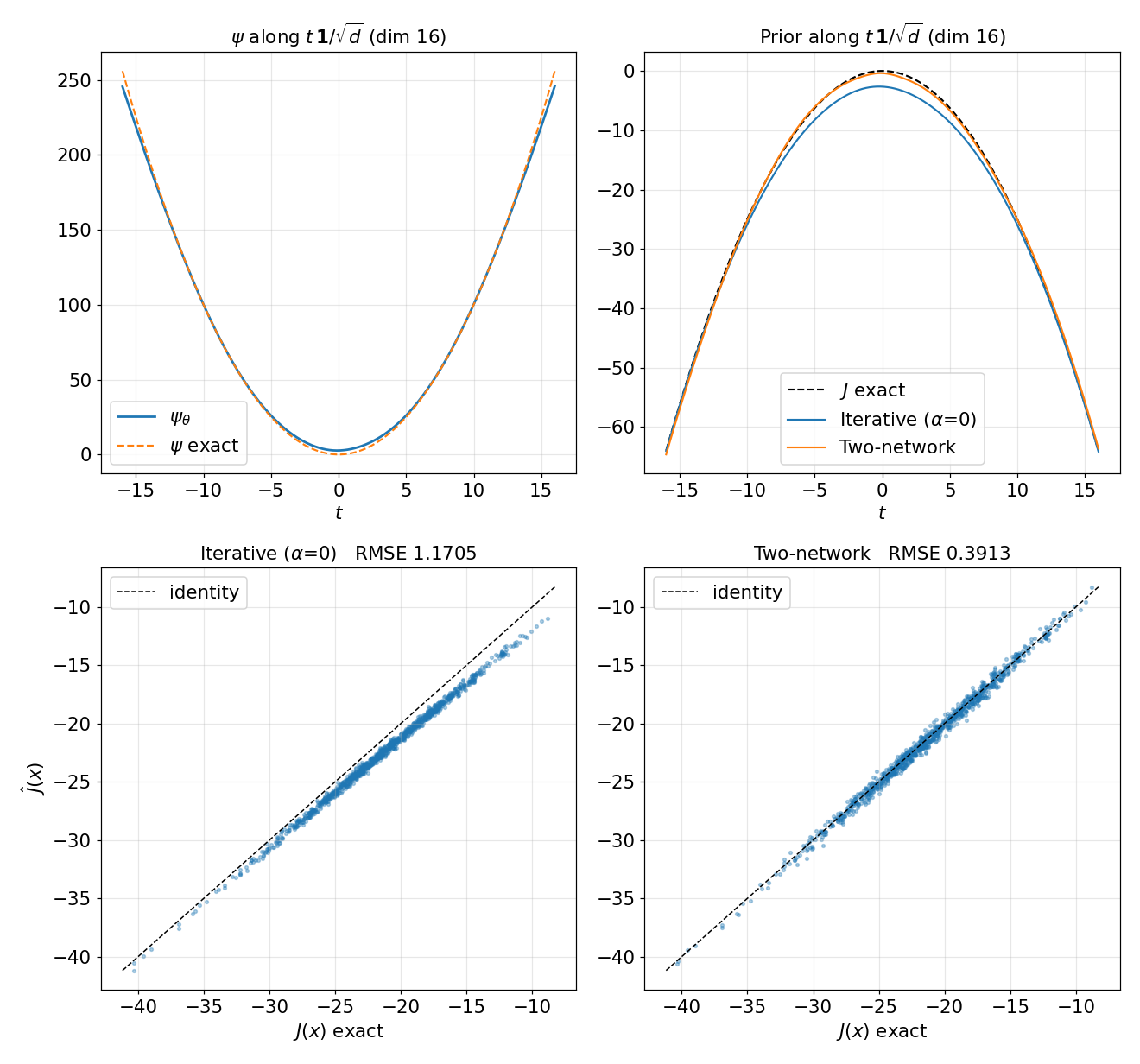}
\caption{\GPLadd{Concave prior, dimension $16$. Panels as in~\cref{dim28}.}}
\label{dim38}
\end{figure}

\GPLadd{\Cref{F3,F3b} report the experiment and~\cref{dim38} plots it at $d = 16$. MAQ is again the most accurate at every dimension, by a factor of $25$ over the best network at $d = 64$.}

\begin{table}[!htbp]
\centering
\caption{\GPLadd{Concave prior: relative $L^{2}$ error of the recovered prior against $\JBVS = J$, columns as in~\cref{F1}. \Cref{ex:quad_exact} applies here with $k = 1$.}}
\footnotesize
\setlength{\tabcolsep}{5pt}
\begin{tabular}{c cccc}
\toprule
$d$ & iterative ($\alpha_{\mathrm{inv}}^{*}$) & two-network & one-network & MAQ\\
\midrule
2 & $5.82\times 10^{-4}$ (0) & $9.30\times 10^{-4}$ & $8.67\times 10^{-4}$ & $\bm{2.43\times 10^{-8}}$\\
4 & $2.07\times 10^{-3}$ (0) & $3.45\times 10^{-3}$ & $3.53\times 10^{-3}$ & $\bm{2.53\times 10^{-7}}$\\
8 & $6.68\times 10^{-3}$ (0) & $9.62\times 10^{-3}$ & $8.96\times 10^{-3}$ & $\bm{3.74\times 10^{-4}}$\\
16 & $4.44\times 10^{-2}$ (0) & $2.23\times 10^{-2}$ & $2.16\times 10^{-2}$ & $\bm{1.39\times 10^{-2}}$\\
32 & $9.48\times 10^{-2}$ (0) & $1.47\times 10^{-1}$ & $1.10\times 10^{-1}$ & $\bm{2.03\times 10^{-2}}$\\
64 & $4.22\times 10^{-1}$ (0.1) & $3.84\times 10^{-1}$ & $3.40\times 10^{-1}$ & $\bm{1.35\times 10^{-2}}$\\
\bottomrule
\end{tabular}
\label{F3}
\end{table}

\begin{table}[!htbp]
\centering
\caption{\GPLadd{Concave prior: validation mean-squared error, as in~\cref{F1b}.}}
\footnotesize
\setlength{\tabcolsep}{4pt}
\begin{tabular}{l cccccc}
\toprule
network & $d = 2$ & $4$ & $8$ & $16$ & $32$ & $64$\\
\midrule
$\psiNN$ & $2.36\times 10^{-6}$ & $7.65\times 10^{-5}$ & $3.46\times 10^{-3}$ & $1.07\times 10^{-1}$ & $2.25\times 10^{0}$ & $4.14\times 10^{1}$\\
$\JNN$ (two-net) & $5.13\times 10^{-6}$ & $2.62\times 10^{-4}$ & $8.87\times 10^{-3}$ & $2.75\times 10^{-1}$ & $5.02\times 10^{0}$ & $3.42\times 10^{1}$\\
$\JNN$ (one-net) & $5.01\times 10^{-6}$ & $2.38\times 10^{-4}$ & $5.41\times 10^{-3}$ & $1.82\times 10^{-1}$ & $4.77\times 10^{0}$ & $3.94\times 10^{1}$\\
\bottomrule
\end{tabular}
\label{F3b}
\end{table}

\subsubsection{Negative $\ell_{1}$ norm}\label{subsubsec:nonnegl1_exp}
\begin{figure}[!htbp]
 \centering
\includegraphics[width=0.8\textwidth,height=0.40\textheight]{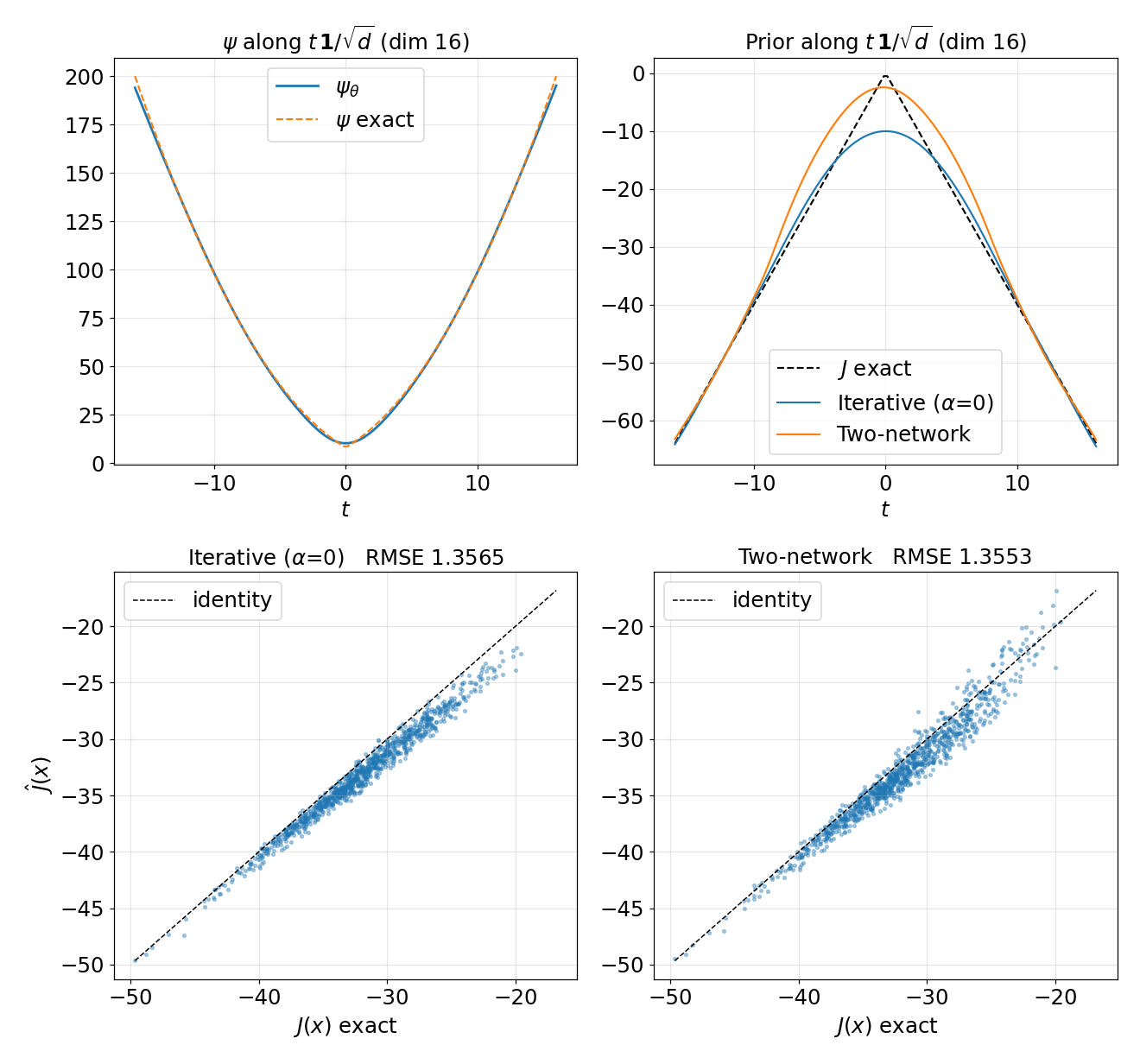}
\caption{\GPLadd{Negative $\ell_1$ prior, dimension $16$. Left: cross sections of the convex potential $\psiNN(\bx;\theta)$ and of the recovered prior, against both references $\JBVS$ and $J$. Right: the root-mean-squared error of the recovered prior against $\JBVS$ on the shared test points, with the iterative recovery of~\cite{fang2024whats} included for comparison.}}
\label{dim841}
\end{figure}

We finish with the nonconvex prior $J(\bx) = -\normone{\bx}$. Since the prior is separable, \cref{ex:negl1_1} gives $S(\bx,t) = -\frac{\GPLadd{d}t}{2} - \normone{\bx}$, and $\JBVS$ is obtained coordinatewise from~\cref{eq:true_jbvs}. \GPLadd{The two references separate most widely on this family, where $\JBVS$ differs from $J$ by between $2.1\%$ and $3.8\%$ in relative $L^{2}$ over the shared test points, the gap narrowing as $d$ grows. \Cref{F4} reports the recovery,~\cref{F4b} the network fits, and~\cref{dim841} plots both at $d = 16$.}

\GPLadd{The iterative inversion is the most accurate at $d = 8$, $16$ and $32$, the two-network recovery at $d = 2$ and $4$, and the one-network variant at $d = 64$, where it beats the two-network recovery by a factor of $6.2$ and the inversion by $3.2$. Two features of the family explain why the network recoveries do worse here than on the other three. First, the map $\nabla_{\bx}\psi(\bx,1) = \bx + \operatorname{sign}(\bx)$ omits the open interval $(-1,1)$ in each coordinate, so the conjugate samples $\by_{k} = \nabla_{\bx}\psiNN(\bx_{k};\theta)$ leave a hole around the origin that no amount of data fills. Second, the prior has a kink at every coordinate hyperplane that a Softplus network smooths. MAQ fails here as well, since the semiconvex transform of $\psi(\cdot,1)$ is a maximum of $2^{d}$ affine pieces.}

\GPLadd{\Cref{F4c} settles which object the data determine, scoring the iterative recovery against each reference. The comparison is informative only where the recovery is more accurate than the separation between the two references, and there it lies between $2.2$ and $4.3$ times closer to $\JBVS$ than to $J$, at every dimension from $2$ to $32$. At $d = 64$ its error of $10\%$ exceeds the $2.1\%$ separation, and the comparison decides nothing.}

\begin{table}[!htbp]
\centering
\caption{\GPLadd{Negative $\ell_{1}$ prior: relative $L^{2}$ error of the recovered prior against $\JBVS$, columns as in~\cref{F1}.}}
\footnotesize
\setlength{\tabcolsep}{5pt}
\begin{tabular}{c cccc}
\toprule
$d$ & iterative ($\alpha_{\mathrm{inv}}^{*}$) & two-network & one-network & MAQ\\
\midrule
2 & $1.12\times 10^{-2}$ (0) & $1.06\times 10^{-2}$ & $\bm{1.04\times 10^{-2}}$ & $1.19\times 10^{-1}$\\
4 & $1.18\times 10^{-2}$ (0) & $\bm{1.15\times 10^{-2}}$ & $1.77\times 10^{-2}$ & $8.85\times 10^{-2}$\\
8 & $\bm{1.36\times 10^{-2}}$ (0) & $1.73\times 10^{-2}$ & $2.04\times 10^{-2}$ & $6.98\times 10^{-2}$\\
16 & $\bm{1.77\times 10^{-2}}$ (0) & $2.81\times 10^{-2}$ & $3.11\times 10^{-2}$ & $4.39\times 10^{-2}$\\
32 & $\bm{1.69\times 10^{-2}}$ (0) & $4.26\times 10^{-2}$ & $2.94\times 10^{-2}$ & $1.81\times 10^{-1}$\\
64 & $1.01\times 10^{-1}$ (0.1) & $1.95\times 10^{-1}$ & $\bm{3.13\times 10^{-2}}$ & $3.20\times 10^{-1}$\\
\bottomrule
\end{tabular}
\label{F4}
\end{table}

\begin{table}[!htbp]
\centering
\caption{\GPLadd{Negative $\ell_{1}$ prior: which reference the data determine. The column $\mathrm{relL}^{2}(J,\JBVS)$ is the separation between the two references on the shared test points, and the next two columns score the iterative recovery against each. A ratio above $1$ means the recovery lies closer to $\JBVS$ than to $J$. We report the ratio only where the recovery error falls below the separation, since a coarser recovery cannot resolve the two; at $d = 64$ the error exceeds the separation.}}
\footnotesize
\setlength{\tabcolsep}{6pt}
\begin{tabular}{c cccc}
\toprule
$d$ & $\mathrm{relL}^{2}(J,\JBVS)$ & iterative vs $\JBVS$ & iterative vs $J$ & ratio\\
\midrule
2 & $3.77\times 10^{-2}$ & $1.12\times 10^{-2}$ & $4.85\times 10^{-2}$ & $4.3$\\
4 & $2.92\times 10^{-2}$ & $1.18\times 10^{-2}$ & $4.08\times 10^{-2}$ & $3.5$\\
8 & $2.65\times 10^{-2}$ & $1.36\times 10^{-2}$ & $4.01\times 10^{-2}$ & $3.0$\\
16 & $2.40\times 10^{-2}$ & $1.77\times 10^{-2}$ & $4.10\times 10^{-2}$ & $2.3$\\
32 & $2.21\times 10^{-2}$ & $1.69\times 10^{-2}$ & $3.63\times 10^{-2}$ & $2.2$\\
64 & $2.12\times 10^{-2}$ & $1.01\times 10^{-1}$ & $8.31\times 10^{-2}$ & ---\\
\bottomrule
\end{tabular}
\label{F4c}
\end{table}

\begin{table}[!htbp]
\centering
\caption{\GPLadd{Negative $\ell_{1}$ prior: validation mean-squared error, as in~\cref{F1b}.}}
\footnotesize
\setlength{\tabcolsep}{4pt}
\begin{tabular}{l cccccc}
\toprule
network & $d = 2$ & $4$ & $8$ & $16$ & $32$ & $64$\\
\midrule
$\psiNN$ & $1.66\times 10^{-4}$ & $1.06\times 10^{-3}$ & $4.71\times 10^{-3}$ & $4.75\times 10^{-2}$ & $5.73\times 10^{-1}$ & $6.84\times 10^{0}$\\
$\JNN$ (two-net) & $2.42\times 10^{-5}$ & $1.67\times 10^{-4}$ & $1.13\times 10^{-3}$ & $5.60\times 10^{-2}$ & $2.21\times 10^{0}$ & $6.36\times 10^{0}$\\
$\JNN$ (one-net) & $2.93\times 10^{-5}$ & $2.04\times 10^{-4}$ & $1.05\times 10^{-3}$ & $4.51\times 10^{-2}$ & $2.83\times 10^{0}$ & $1.59\times 10^{1}$\\
\bottomrule
\end{tabular}
\label{F4b}
\end{table}

\subsubsection{\GPLadd{Comparison across the four families}}\label{subsubsec:cross_family}
\begin{figure}[!htbp]
\centering
\includegraphics[width=\textwidth,height=0.40\textheight,trim=0 0 0 70,clip]{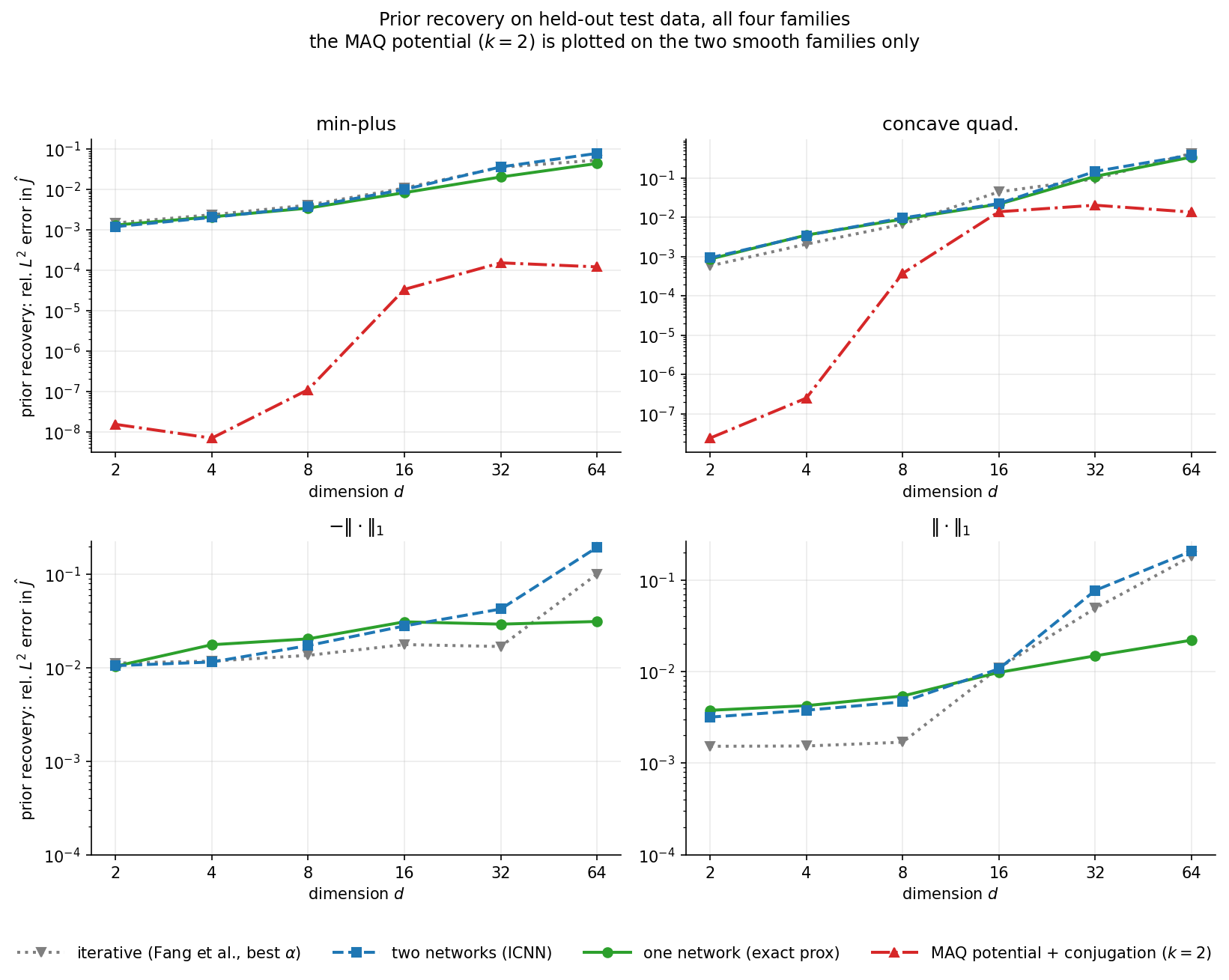}
\caption{\GPLadd{Prior recovery against dimension for the four families of~\cref{subsec:closed_form}, measured as the relative $L^{2}$ error~\cref{eq:relL2} against $\JBVS$ on the $1000$ shared test points. The MAQ potential is plotted on the two families of~\cref{ex:quad_exact} only.}}
\label{fig:summary_recovery}
\end{figure}

\GPLadd{\Cref{fig:summary_recovery} plots the recovery error against the dimension for the four families of~\cref{subsec:closed_form}. Every panel carries the three network recoveries, and the two panels covered by~\cref{ex:quad_exact} also carry MAQ, the most accurate method there and the least accurate on the other two families. The one-network and two-network recoveries agree to within a factor of about $1.5$ through $d = 16$ on every family, and they separate above it. The one-network error grows more slowly in $d$ than the two-network error on every family, and the two separate above $d = 16$. In~\cref{fig:ray_minplus} the first network loses accuracy as the dimension grows, and the second network, which trains on the samples the first one produces, inherits that error.}

\GPLadd{The iterative recovery of~\cite{fang2024whats} is the most accurate method at the smaller dimensions on three of the four families. It degrades at $d = 64$, where the best regularization parameter switches from $\alpha_{\mathrm{inv}} = 0$ to $\alpha_{\mathrm{inv}} = 0.1$ on all four families at once.}

\GPLadd{The references $J$ and $\JBVS$ differ by a few percent on $-\normone{\bx}$, and~\cref{F4c} settles the question there. On the min-plus prior the separation runs from $4.7\times10^{-4}$ at $d = 2$ to $1.2\times10^{-4}$ at $d = 64$, and no network recovery resolves it. The MAQ recovery does. At $d = 16$ it lies $3.4\times10^{-5}$ from $\JBVS$ and $3.0\times10^{-4}$ from $J$, an order of magnitude below the separation. Every recovery accurate enough to tell the two references apart agrees with $\JBVS$, as~\cref{thm:prior_optimality} predicts.}

\begin{figure}[!htbp]
\centering
\includegraphics[width=\textwidth,height=0.40\textheight,trim=0 0 0 70,clip]{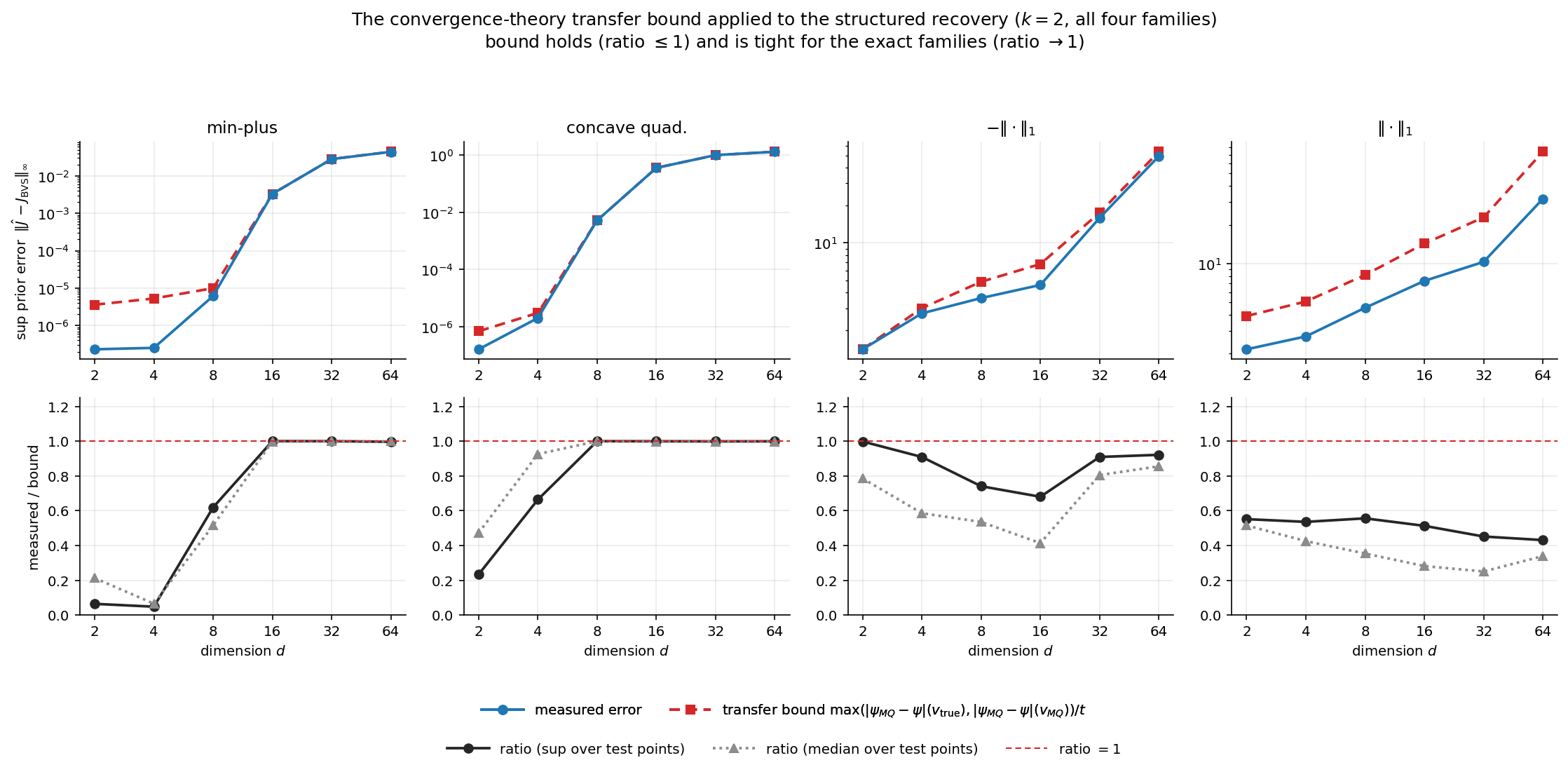}
\caption{\GPLadd{Measured prior-recovery error of the MAQ recovery against the pointwise transfer bound~\cref{eq:pointwise_bound}, for each family and dimension.}}
\label{fig:theory_bounds}
\end{figure}

\GPLadd{\Cref{fig:theory_bounds} tests the transfer bound of~\cref{prop:transfer_values} against measurement. We apply it to MAQ, whose potential error $|\psi_{\mathrm{MAQ}} - \psi(\cdot,t)|$ we can evaluate at the two preimages~\cref{prop:transfer_values} identifies, so we compare against the pointwise bound~\cref{eq:pointwise_bound} itself rather than a numerically estimated supremum. On the min-plus and concave priors, the ratio of measured error to bound agrees with $1$ to two decimals from $d = 16$ upward. On the two kink priors, whose potentials $k = 2$ atoms cannot represent, the ratio stays between $0.43$ and $1.00$.}

\subsection{\GPLadd{The implicit prior of a total-variation posterior-mean denoiser}}\label{subsec:tv_exp}
\GPLadd{Darbon and Langlois~\cite{darbon2021bayesian} showed that a total variation posterior mean denoiser is the proximal mapping of a twice continuously differentiable convex prior. We learn that prior from data. Write $\bu \in \R^{n}$ for an image patch with $n$ pixels arranged on a square grid, and let
\begin{equation}\label{eq:atv}
\TV(\bu) \coloneqq \sum_{(i,j) \in \mathcal{N}} \left|u_{i} - u_{j}\right|
\end{equation}
denote anisotropic total variation, the sum over unordered horizontally or vertically adjacent pixel pairs. Given noisy data $\bx$ and a denoising strength $t > 0$, set
\[
E(\bu;\bx) \coloneqq \frac{1}{2t}\normsq{\bu-\bx} + \TV(\bu),
\]
and let the \emph{Gibbs posterior} at \emph{temperature} $\epsilon > 0$ be the density on $[0,1]^{n}$ proportional to $e^{-E(\cdot\,;\bx)/\epsilon}$. Its mode is the maximum a posteriori (MAP) estimate $\prox_{t\TV}(\bx)$, and we study its \emph{posterior mean}
\begin{equation}\label{eq:upm}
\uPM(\bx) \coloneqq \frac{\int_{[0,1]^{n}} \bu\, e^{-E(\bu;\bx)/\epsilon}\diff\bu}{\int_{[0,1]^{n}} e^{-E(\bu;\bx)/\epsilon}\diff\bu}.
\end{equation}
}

\subsubsection{\GPLadd{The posterior mean as a proximal operator}}\label{subsubsec:tv_theory}
\GPLadd{Define
\begin{equation}\label{eq:cole_hopf}
\Seps(\bx,t) \coloneqq -\epsilon\ln\left(\frac{1}{(2\pi t \epsilon)^{n/2}}\int_{\Rn} e^{-\left(\frac{1}{2t}\normsq{\bx-\by} + J(\by)\right)/\epsilon}\diff\by\right),
\end{equation}
with $J = \TV$ plus the indicator of $[0,1]^{n}$, the effective prior of~\cref{eq:upm}. The Cole--Hopf transformation turns the linear heat equation satisfied by the integral in~\cref{eq:cole_hopf} into the \emph{viscous} Hamilton--Jacobi equation
\begin{equation}\label{eq:viscous_hj_tv}
\frac{\partial \Seps}{\partial t}(\bx,t) + \frac{1}{2}\normsq{\nabla_{\bx}\Seps(\bx,t)} = \frac{\epsilon}{2}\Delta_{\bx}\Seps(\bx,t), \qquad \Seps(\bx,0) = J(\bx),
\end{equation}
which $\Seps$ solves uniquely and smoothly~\cite[\S3]{darbon2021bayesian}. It is the viscous counterpart of~\cref{eqn:intro2}, and differentiating~\cref{eq:cole_hopf} under the integral sign identifies the posterior mean as
\begin{equation}\label{eq:pm_is_grad}
\uPM(\bx) = \bx - t\nabla_{\bx}\Seps(\bx,t) = \nabla_{\bx}\Keps(\bx,t), \qquad \Keps(\bx,t) \coloneqq \frac{1}{2}\normsq{\bx} - t\Seps(\bx,t),
\end{equation}
with $\bx \mapsto \Keps(\bx,t)$ strictly convex~\cite[Thm.~3.3]{darbon2021bayesian}. The map $\bx \mapsto \Seps(\bx,t)$ is convex~\cite[\S3]{darbon2021bayesian}, which gives $\nabla_{\bx}\uPM = \matr{I}_{n \times n} - t\nabla_{\bx}^{2}\Seps(\bx,t) \preceq \matr{I}_{n \times n}$, while the strict convexity of $\bx \mapsto \Keps(\bx,t)$ quoted above supplies the matching lower bound $\nabla_{\bx}\uPM = \nabla_{\bx}^{2}\Keps(\bx,t) \succeq 0$; together $0 \preceq \nabla_{\bx}\uPM \preceq \matr{I}_{n \times n}$, so $\uPM$ is nonexpansive and the second part of~\cref{thm:char_prox} with $L = 1$ makes the recovered prior convex. Hence $\uPM$ is the proximal operator of
\begin{equation}\label{eq:freg}
\freg \coloneqq \Keps^{*}(\cdot,t) - \frac{1}{2}\normsq{\cdot},
\end{equation}
and $\freg$ is twice continuously differentiable even though $\TV$ is not.}

\GPLadd{The substitutions
\[
S(\bx,t) \longleftarrow \Seps(\bx,t), \qquad \psi(\bx,t) \longleftarrow \Keps(\bx,t), \qquad t\JBVS \longleftarrow \freg
\]
carry~\cref{eq:pm_is_grad,eq:freg} onto~\cref{eq:psi_defn,eq:formula_recovery_prior}, so the viscous problem folds onto the first-order problem this paper already treats and the recovery identity~\cref{eq:formula_recovery_prior} becomes~\cref{eq:freg}. Because $\freg$ is smooth, the posterior-mean estimator avoids \emph{staircasing}, the artifact by which a total-variation MAP estimate renders a smooth gradient as flat steps separated by jumps.}

\GPLadd{In what follows, we describe the training procedure for learning the implicit prior associated to total variation posterior mean estimation.}

\GPLadd{\emph{Parameters.} The noise level $\sigma$ and the denoising strength $t$ determine the temperature through $\epsilon = \sigma^{2}/t$. We report $\sigma = t = 20/256 \approx 0.0781$ throughout, in units where pixel values lie in $[0,1]$, and work with $8\times8$ patches, so $n = 64$. The box $[0,1]^{n}$ is part of the model, so the effective prior is $\TV$ plus the indicator of $[0,1]^{n}$, the denoiser $\uPM$ maps into the open box, and $\freg$ is defined there.}

\GPLadd{\emph{Evaluating the denoiser.} Since~\cref{eq:upm} has no closed form at $n = 64$, we evaluate $\uPM$ by Markov chain Monte Carlo, using single-site random-walk Metropolis--Hastings with a uniform proposal and pixels clipped to $[0,1]$, one independent chain per patch, started at $\bu^{0} = \bx$ and advanced by $m = 8000$ \emph{sweeps}, where one sweep is one proposed update at every pixel. Because $E(\cdot\,;\bx)$ couples only adjacent pixels, the pixels of one color of the checkerboard are conditionally independent given the other, so the sweep updates one color at a time across all chains at once, which changes the scan order and nothing else. We take $\uPM(\bx)$ to be the running mean of the chain after discarding its \emph{burn-in}, the early portion in which the chain still remembers $\bu^{0}$ rather than sampling from the posterior.}

\GPLadd{\emph{Data.} We draw $20\,000$ training patches at random positions from one image (Barbara), add independent Gaussian noise of standard deviation $\sigma$, clip to $[0,1]$, and evaluate $\uPM$ at each, with seed $1$. The validation set is $4000$ patches at positions distinct from the training ones in the same image, with seed $2$, and the evaluation set is $4000$ patches from a \emph{different} image, cameraman, with seed $3$, of which the first $2000$ are scored. Training and validation patches overlap in pixels, since a $256\times256$ image holds only $1024$ disjoint $8\times8$ patches, so validation measures generalization over the noise and the sampler while the evaluation split carries the test over image content. Every number below is therefore a \emph{transfer} number, computed on an image the networks never saw, since $\freg$ is a property of the denoiser and not of an image.}

\GPLadd{\emph{Training.} The training target is the proximal optimality condition
\begin{equation}\label{eq:prox_opt_tv}
\nabla\freg(\by_{k}) = \bx_{k} - \by_{k}, \qquad \by_{k} \coloneqq \uPM(\bx_{k}).
\end{equation}
The loss is the mean-squared error of $\nabla J_{\theta}(\by_{k})$, the gradient of the one-network prior $J_{\theta}$, against that right-hand side, normalized by $\operatorname{var}(\bx-\by)$. The architecture, optimizer, batch size, budget and schedule are those of~\cref{subsubsec:protocol}, with two changes. First, we set the Softplus parameter to $\beta = 20$ rather than the $5$ used in~\cref{subsec:closed_form}, since $\freg$ inherits from $\TV$ a kink rounded over a scale we measure at $n = 2$ to be about $0.038$, which $\beta = 5$ is too smooth to resolve. Second, we standardize the input, because patch values concentrate on a thin set (standard deviation $0.039$ within a patch) while Softplus bends over a scale of about $0.2$, leaving the unstandardized network nearly affine in the informative directions.}

\GPLadd{\emph{The score.} With no analytic prior to compare against, we score a recovered prior $\hat{J}$ by its \emph{held-out proximal residual}, the relative failure of $\nabla\hat{J}$ to satisfy~\cref{eq:prox_opt_tv} at an evaluation patch,
\begin{equation}\label{eq:tv_resid}
r(\bx_{k}) \coloneqq \frac{\norm{\nabla\hat{J}(\by_{k}) - (\bx_{k}-\by_{k})}}{\norm{\bx_{k}-\by_{k}}}, \qquad \by_{k} = \uPM(\bx_{k}).
\end{equation}
The iterative recovery supplies $\nabla\hat J(\by) = \hat{\bx}(\by) - \by$, where $\hat{\bx}(\by)$ is the preimage solving $\nabla\psiNN(\hat{\bx}) = \by$ and is the point $\hat{\bx}$ of~\cref{prop:transfer_values}. The two-network recovery supplies $\nabla\JNN(\by;\tilde{\theta}) - \by$, and the one-network recovery supplies $\nabla J_{\theta}(\by)$. Because $\uPM$ is itself estimated by a finite chain, $r$ has a floor. We put the sampler's own relative noise at
\[
\delta \coloneqq \frac{1}{16}\sum_{k=1}^{16} \frac{\bigl\|\operatorname{std}_{1 \leqslant j \leqslant 8}\bigl(\bx_{k}-\by_{k}^{(j)}\bigr)\bigr\|}{\norm{\bx_{k}-\by_{k}}} = 2.19\%,
\]
where $\operatorname{std}$ is taken coordinate-wise over $8$ independent chains per patch and the average runs over $16$ evaluation patches. No method should fall below $\delta$.}

\subsubsection{\GPLadd{TV Posterior mean estimator: Results}}\label{subsubsec:tv_results}

\begin{table}[!htbp]
\centering
\caption{\GPLadd{Recovery of the implicit regularizer $\freg$ of the total-variation posterior-mean denoiser. Each entry is the held-out proximal residual~\cref{eq:tv_resid}, over $2000$ evaluation patches of a transfer image; ``median'' and ``p90'' are its median and ninetieth percentile, and ``$\times\delta$'' is the median divided by the sampler noise floor $\delta = 2.19\%$, the relative spread of the target across independent chains. ``Cost'' is the wall-clock time to evaluate the prior at all $2000$ patches on one CPU. The iterative recovery is run unregularized, at $\alpha_{\mathrm{inv}} = 0$. Boldface marks the best entry in a column.}}
\footnotesize
\begin{tabular}{l cccc}
\toprule
\GPLadd{Recovered prior} & \GPLadd{Residual, median} & \GPLadd{Residual, p90} & \GPLadd{$\times\delta$} & \GPLadd{Cost}\\
\midrule
\GPLadd{Iterative~\cite{fang2024whats}} & \GPLadd{$11.8\%$} & \GPLadd{$377\%$} & \GPLadd{$5.4$} & \GPLadd{$3.2$ s (one inversion per patch)}\\
\GPLadd{Two-network} & \GPLadd{$13.7\%$} & \GPLadd{$98\%$} & \GPLadd{$6.3$} & \GPLadd{$0.01$ s (forward pass)}\\
\GPLadd{One-network} & \GPLadd{$\bm{8.1\%}$} & \GPLadd{$\bm{21\%}$} & \GPLadd{$\bm{3.7}$} & \GPLadd{$0.03$ s (forward pass)}\\
\bottomrule
\end{tabular}
\label{tab:tv_recovery}
\end{table}

\begin{figure}[!htbp]
\centering
\makebox[\textwidth][c]{\includegraphics[width=1\textwidth]{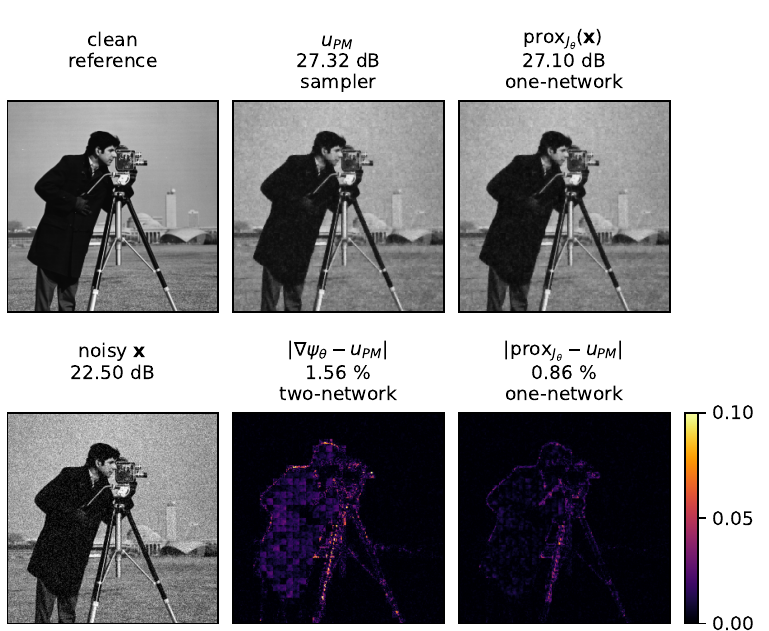}}
\caption{\GPLadd{The proximal map of the recovered prior, on the transfer image of~\cref{tab:tv_denoise}. Top: the clean image, the posterior-mean denoiser $\uPM(\bx)$ produced by the sampler, and the output $\prox_{J_{\theta}}(\bx)$ of the proximal map of the one-network recovered prior, scored in decibels against the clean image. Bottom: the noisy input $\bx$, and the errors of $\nabla\psiNN$ and of $\prox_{J_{\theta}}$ against $\uPM$ on a common magnified scale, in relative $L^{2}$. Both denoised images are indistinguishable from $\uPM$ at this scale, and both errors concentrate at edges, where $\freg$ is most nonlinear.}}
\label{fig:tv_denoise}
\end{figure}

\GPLadd{The one-network recovery is the most accurate of the three, at $3.7$ times the sampler floor against $6.3$ for the two-network recovery and $5.4$ for the iterative recovery, and it has the smallest tail, with a ninetieth percentile of $21\%$ against $98\%$ and $377\%$.\footnote{\GPLadd{Retraining $\psiNN$ and $\JNN$ and repeating the comparison leaves the one-network column unchanged at $8.1\%$, $21\%$ and $3.7$, and moves the two baselines. The iterative recovery then reads $12.9\%$, $60\%$ and $5.9$, and the two-network recovery $13.4\%$, $64\%$ and $6.1$. In~\cref{tab:tv_denoise} the same retraining moves $1.56\%$ to $1.39\%$ and leaves $0.86\%$ unchanged. The inversion diverges in both runs. This table and~\cref{tab:tv_denoise,fig:tv_denoise} report the same run.}} Its values correlate with $\TV$ over the evaluation patches at a Pearson coefficient of $0.983$. We run the iterative recovery at $\alpha_{\mathrm{inv}} = 0$, where the objective need not be coercive, and on the patches where it fails the iterate leaves the region in which the trained potential is accurate, returning a largest preimage coordinate of $2.27$ in absolute value against the range $[0,1]$ in which the patches live.}

\GPLadd{The residual tests only a first-order condition, so we also check the recovery directly. If we have recovered $\freg$ faithfully, the proximal map of the learned prior must reproduce the denoiser.}

\begin{table}[!htbp]
\centering
\caption{\GPLadd{Agreement between the proximal map of the recovered prior and the posterior-mean denoiser $\uPM$, on a full transfer image (cameraman, $256\times256$, $\sigma = t = 20/256$, $m = 8000$ sweeps, seed $3$). The three error columns take $\uPM$ as the reference, not the clean image: ``rel.\ $L^{2}$'' is $\norm{\hat\bu - \uPM}/\norm{\uPM}$, PSNR is the peak signal-to-noise ratio $10\log_{10}(1/\mathrm{MSE})$, reported in decibels (dB) with pixel values in $[0,1]$, and SSIM is the structural similarity index, which equals $1$ for identical images. The reference block scores the same maps against the clean image instead. Boldface marks the better of the two recoveries.}}
\footnotesize
\begin{tabular}{l ccc cc}
\toprule
& \multicolumn{3}{c}{\GPLadd{Against $\uPM$}} & \multicolumn{2}{c}{\GPLadd{Against the clean image}}\\
\cmidrule(lr){2-4}\cmidrule(lr){5-6}
\GPLadd{Map} & \GPLadd{rel.\ $L^{2}$} & \GPLadd{PSNR} & \GPLadd{SSIM} & \GPLadd{PSNR} & \GPLadd{SSIM}\\
\midrule
\GPLadd{Noisy input $\bx$} & \GPLadd{---} & \GPLadd{---} & \GPLadd{---} & \GPLadd{$22.50$ dB} & \GPLadd{$0.411$}\\
\GPLadd{$\uPM$, the sampler} & \GPLadd{---} & \GPLadd{---} & \GPLadd{---} & \GPLadd{$27.32$ dB} & \GPLadd{$0.715$}\\
\midrule
\GPLadd{$\nabla\psiNN$ (iterative and two-network)} & \GPLadd{$1.56\%$} & \GPLadd{$41.76$ dB} & \GPLadd{$0.9945$} & \GPLadd{$26.92$ dB} & \GPLadd{$0.704$}\\
\GPLadd{One-network, $\prox_{J_{\theta}}$} & \GPLadd{$\bm{0.86\%}$} & \GPLadd{$\bm{46.93}$ dB} & \GPLadd{$\bm{0.9982}$} & \GPLadd{$\bm{27.09}$ dB} & \GPLadd{$\bm{0.711}$}\\
\bottomrule
\end{tabular}
\label{tab:tv_denoise}
\end{table}

\GPLadd{Note that neither learned-potential recovery inverts anything to \emph{denoise}. Since $\nabla\psi^{*}(\cdot,t) = (\nabla_{\bx}\psi(\cdot,t))^{-1}$, both deploy the same map $\nabla\psiNN$ and share a row, differing only in how they \emph{evaluate} the prior, which~\cref{tab:tv_recovery} measures. The one-network recovery reproduces the denoiser to $0.86\%$ in relative $L^{2}$ against $1.56\%$, and both fall far below the $10\%$ scale of the residuals in~\cref{tab:tv_recovery}. There is no contradiction. Those residuals measure agreement of the \emph{gradients} of the priors; this table measures agreement of the resulting maps, and the proximal map contracts the discrepancy.}

\subsection{Plug-and-play priors}\label{subsec:pnp_exp}
Most iterative solvers for inverse problems involve the regularizer only through its proximal operator, and by~\cref{eq:def_prox_argmin} that step is the maximum a posteriori (MAP) estimate for $\by$ observed in Gaussian noise of variance $t$. Replacing it by an off-the-shelf denoiser, so that no $J$ is ever formed, is the plug-and-play (PnP) framework of Venkatakrishnan, Bouman and Wohlberg~\cite{venkatakrishnan2013plug}, developed with learned denoisers in~\cite{zhang2017learning, meinhardtlearning, zhang2021plug} and surveyed in~\cite{kamilov-2023-plug}.

PIRATE~\cite{hu2023plug} does this for deformable image registration with a learned operator $\mathsf{D}$, minimizing a data term $g$ plus a regularizer $R$ and updating the registration field $\phi$ by
\begin{equation}\label{eq:pirate-update}
 \phi \;\leftarrow\; \phi - \gamma\bigl(\nabla g + \alpha_{\mathrm{R}}\nabla R
 + \tau\,(\phi - \mathsf{D}(\phi))\bigr),
\end{equation}
with step size $\gamma$, regularization weight $\alpha_{\mathrm{R}}$ and denoiser weight $\tau$ so the implicit regularizer is the term whose gradient is the residual $\phi-\mathsf{D}(\phi)$. That the residual is a gradient is justified through the identity for the exact minimum-mean-squared-error (MMSE) denoiser, but the operator deployed is a trained network, and whether its residual is the gradient of anything is neither stated nor tested.

\Cref{thm:char_prox} gives $\mathsf{D}=\prox_{tJ}$ for some proper $J$ only if $\mathsf{D}(\bx)\in\partial_{\bx}\psi(\bx,t)$ for a convex $\psi$. 
The Jacobian $\nabla \mathsf{D}(\bx)$ therefore must be symmetric and positive semidefinite, with $\nabla \mathsf{D}\preceq \matr{I}_{n \times n}$ exactly when the regularizer is convex.

Now, writing $\by_k=\mathsf{D}(\bx_k)$ and using~\cref{eqn:early_prox}, each denoiser evaluation supplies a gradient sample
\begin{equation}\label{eq:pnp-grad}
 \nabla_{\bx}S(\bx_k,t)=\frac{\bx_k-\by_k}{t}
 \quad \iff\quad
 \nabla_{\by}\bigl(t\JBVS\bigr)(\by_k)=\bx_k-\by_k.
\end{equation}
So the PnP setting is the incomplete-information problem of~\cref{sec:inverse_incomplete_info} with gradient data in place of values. By~\cref{thm:prior_optimality} what~\cref{eq:pnp-grad} determines is $\JBVS$, and $J\geqslant\JBVS$ in general, so a denoiser reveals $\freg\coloneqq t\JBVS$ and not $tJ$.

The noise level is a time, and the viscous problem of~\cref{subsubsec:tv_theory} returns here with its temperature fixed at $\epsilon = 1$. Let the clean field have prior density $p_{\phi}\propto e^{-J}$, let $\bx=\phi+\sigma\bm{\xi}$ with $\bm{\xi}\sim\mathcal{N}(0,\matr{I}_{n \times n})$, and write $\Dstar(\bx)\coloneqq\mathbb{E}[\phi\mid \bx]$ for the exact MMSE denoiser. Set $t=\sigma^{2}$ so that $\epsilon=\sigma^{2}/t$ equals $1$, and write $G_{t}$ for the heat kernel at time $t$. The noisy density is then $p_{\bx}=p_{\phi} * G_{t}$, and the value function~\cref{eq:cole_hopf} built from this $J$ satisfies $e^{-\Seps(\cdot,t)}=G_{t} * e^{-J}\propto p_{\bx}$, so~\cref{eq:viscous_hj_tv} reads
\begin{equation}\label{eq:viscous-hj}
 \frac{\partial \Seps}{\partial t}(\bx,t)+\frac{1}{2}\normsq{\nabla_{\bx}\Seps(\bx,t)} = \frac{1}{2}\Delta_{\bx}\Seps(\bx,t),
 \qquad \Seps(\bx,0)=J(\bx).
\end{equation}
Moreover, note that
\begin{equation}\label{eq:tweedie}
 \Dstar(\bx)=\bx-t\nabla_{\bx}\Seps(\bx,t)=\nabla_{\bx}\Keps(\bx,t),
 \qquad \Keps(\bx,t) \coloneqq \frac{1}{2}\normsq{\bx}-t\Seps(\bx,t).
\end{equation}
Thus $\Dstar$ is the time-$t$ flow map of~\cref{eq:viscous-hj} and its residual $\bx-\Dstar(\bx)=t\nabla_{\bx}\Seps(\bx,t)$ is a gradient. The covariance identity $\nabla \Dstar=\nabla_{\bx}^{2}\Keps(\bx,t)=\operatorname{Cov}(\phi\mid \bx)/\sigma^{2}\succeq0$ makes $\bx \mapsto \Keps(\bx,t)$ convex for \emph{every} prior, with no log-concavity assumed. So~\cref{thm:char_prox} gives $\Dstar=\prox_{\freg}$ with $\freg=\Keps^{*}(\cdot,t)-\tfrac12\normsq{\cdot}$, which is~\cref{eq:freg} and agrees with $\freg=t\JBVS$ by~\cref{eq:formula_recovery_prior}.

The parallel with~\cref{subsubsec:tv_theory} stops at convexity. There $\bx \mapsto \Seps(\bx,t)$ was itself convex, which forced $\nabla\uPM\preceq\matr{I}_{n \times n}$ and left $\freg$ convex; a nonconvex $J$ carries no such guarantee, and all that survives is $\nabla^{2}\freg=\nabla^{2}\Keps^{*}(\cdot,t)-\matr{I}_{n \times n}\succeq-\matr{I}_{n \times n}$. The recovered $\freg$ is therefore $1$-semiconvex and $\JBVS=\freg/t$ has curvature $\succeq-1/t$.

We test the recovery on a Gaussian mixture prior in one dimension, whose posterior-mean denoiser $\mathsf{D}$, value function $S$, potential $\psi$ and backward solution $\JBVS$ are all available in closed form. We draw $N=2\times10^{4}$ samples at $\sigma=0.5$ with $t=\sigma^{2}$. Every network here is supervised on gradients as in~\cref{subsec:tv_exp}, not on values as in~\cref{subsubsec:two_phase}, and training otherwise follows~\cref{subsubsec:protocol} at $4\times10^{4}$ steps with $2$ layers, $64$ hidden units and Softplus parameter $\beta=20$.

\begin{table}[!htbp]
 \centering
 \caption{\GPLadd{One-dimensional mixture denoiser, $\sigma=0.5$, all networks $4\times10^{4}$ steps. Errors are relative $L^{2}$ against the exact $\freg$ on the sampled range of $\mathsf{D}$, each side centered, since gradient supervision does not identify the additive constant; the prox residual is the median of $\|\nabla f(\by_k)-(\bx_k-\by_k)\|/\|\bx_k-\by_k\|$ on held-out data. The exact minimum curvature is $-0.9113$, against the semiconvexity floor $-1$.}}
 \label{tab:gaussian_mixture}
 \begin{tabular}{llrrrr}
 \hline
 recovery & class & rel.\ $L^{2}$ & min.\ curv. & prox resid. & eval (s)\\
 \hline
 Two-network & $1$-semiconvex & $0.17\%$ & $-0.9619$ & $0.08\%$ & $0.029$\\
 Iterative & $1$-semiconvex & $0.03\%$ & $-0.9120$ & $0.17\%$ & $3.251$\\
 One-network & convex & $101.21\%$ & $-0.0048$ & $59.37\%$ & $0.028$\\
 \hline
 \end{tabular}
\end{table}

\begin{figure}[!htbp]
 \centering
\includegraphics[width=1\textwidth,height=0.25\textheight,trim=0 0 0 90,clip]{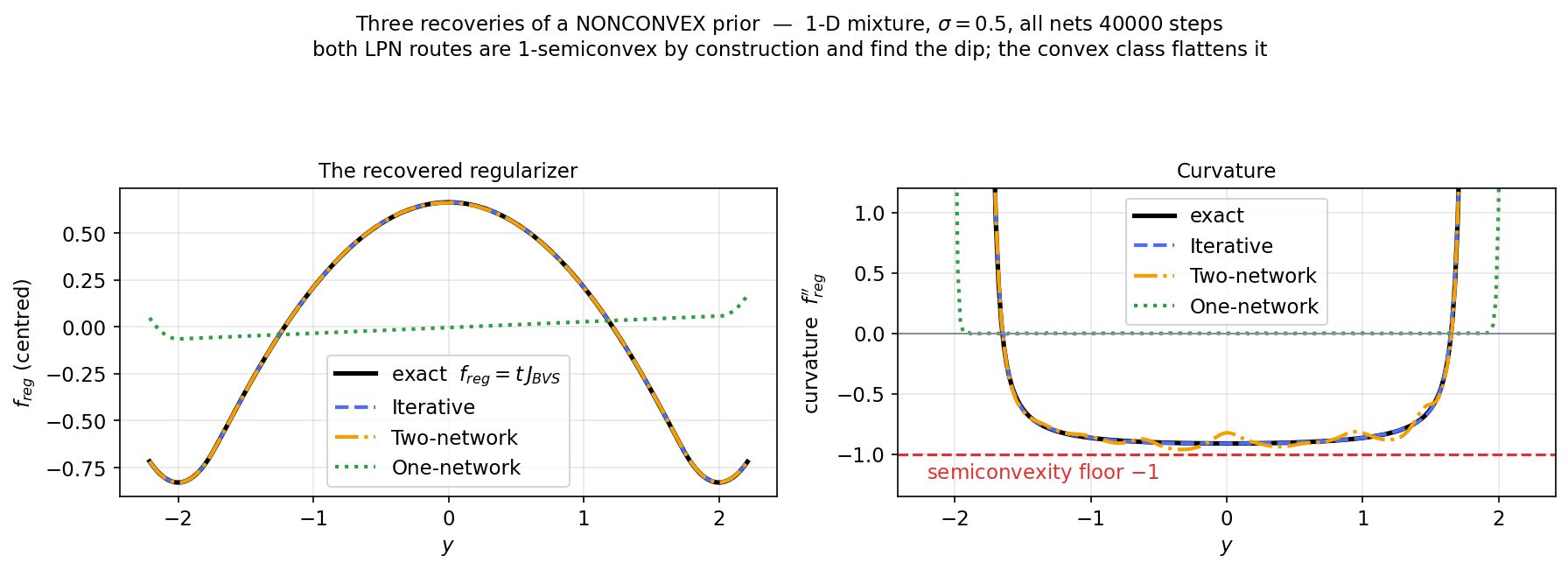}
 \caption{\Toadd{Prior of a one-dimensional Gaussian mixture, $\freg=t\JBVS$ (black), $\sigma=0.5$. \emph{Left}: recovered regularizer. \emph{Right}: curvature, against the semiconvexity floor $-1$.}}
 \label{fig:gaussian_mixture}
\end{figure}

In~\cref{tab:gaussian_mixture}, the one-network recovery reads a minimum
curvature of $-0.0048$ where the exact value is $-0.9113$. Its best validation
loss is $4.8\times10^{-1}$, against $6.5\times10^{-9}$ for the iterative
recovery and $2.1\times10^{-7}$ for the two-network recovery under an identical
budget. The optimizer is therefore not the limitation; the function class is. 
The two semiconvex recoveries agree with the exact prior and with each other to
within $0.14$ percentage points. The inversion in iterative recovery is exact to machine precision, so the
agreement is a statement about the recovery mathematics and not about a solver.
What separates them is the cost. The forward pass evaluates the recovery in $0.029$~seconds against
$3.251$~seconds by the iterative method, which is about $110$ times faster.

\Cref{fig:gaussian_mixture} shows the same three curves. In the left panel the
two semiconvex recoveries are indistinguishable from the exact $\freg$ across the
sampled range, while the one-network recovery is flat. The one-network recovery fits a plain ICNN,
whereas the recoverable object $\freg=\psi^{*}-q$, with
$q=\tfrac12\|\cdot\|_2^2$, satisfies $\nabla^{2}\freg=\nabla^{2}\psi^{*}-I
\succeq-I$ and is only $1$-semiconvex. The hypothesis class therefore does not contain the target, and fitting
returns its convex projection. It
reproduces neither the bump between the modes nor the wells on either side. On the
right panel, the exact curvature falls to a plateau near
$-0.91$ in the interior and rises steeply outside it. Both semiconvex curves
follow it down, but the one-network curve stays at zero. 

We repeat the comparison at $n=64$ on a bimodal field prior on $8\times8$ patches, nonconvex along a single bump-pattern direction $\bu$, along which it reduces to the mixture above, and Gaussian with spectrum $\lambda_k=4/k^{2}$ on the remaining $63$ directions. The backward solution is known in closed form, hull-limited along $\bu$ and exactly $tJ$ on the convex directions. Here the two-network recovery trains $\JNN$ on the denoiser's own outputs $\by_{k}=\mathsf{D}(\bx_{k})$ rather than on $\nabla\psiNN(\bx_{k})$ as in~\cref{subsubsec:protocol}, since at $n=64$ the conjugate carries $\psiNN$'s error into a systematic offset in the recovered values, at a cost of a factor of $19$. Training is as above on $N=4\times10^{4}$ samples at $5\times10^{4}$ steps with $256$ hidden units and standardized inputs, and the inversion is the Adam solver of~\cite{fang2024whats} run unregularized, at $\alpha_{\mathrm{inv}}=0$, for $2\times10^{4}$ iterations per query block.

\begin{table}[!htbp]
 \centering
 \caption{\Toadd{Bimodal field prior, $n=64$, $\sigma=0.5$. ``Value'' is relative $L^{2}$ against the exact backward solution, centered as in~\cref{tab:gaussian_mixture}; ``$\bu$ curv.'' is the minimum curvature along the nonconvex direction (exact: $-0.9113$); ``$\perp$ spec.'' is the relative error of the $63$ Gaussian-mode curvatures; ``resid.'' is the median prox residual on held-out data.}}
 \label{tab:pnp-bimodal}
 \begin{tabular}{llrrrr}
 \hline
 recovery & class & value & $\bu$ curv. & $\perp$ spec. & resid.\\
 \hline
 Two-network & $1$-semiconvex & $3.01\%$ & $-0.9298$ & $0.43\%$ & $3.21\%$\\
 Iterative & $1$-semiconvex & $1.56\%$ & $-0.8600$ & $14.62\%$ & $12.97\%$ \\
 One-network & convex & $15.92\%$ & $+0.0010$ & $1.27\%$ & $6.00\%$\\
 \hline
 \end{tabular}
\end{table}

\begin{figure}[!htbp]
 \centering
 \includegraphics[width=1\textwidth,height=0.25\textheight,trim=0 0 0 70,clip]{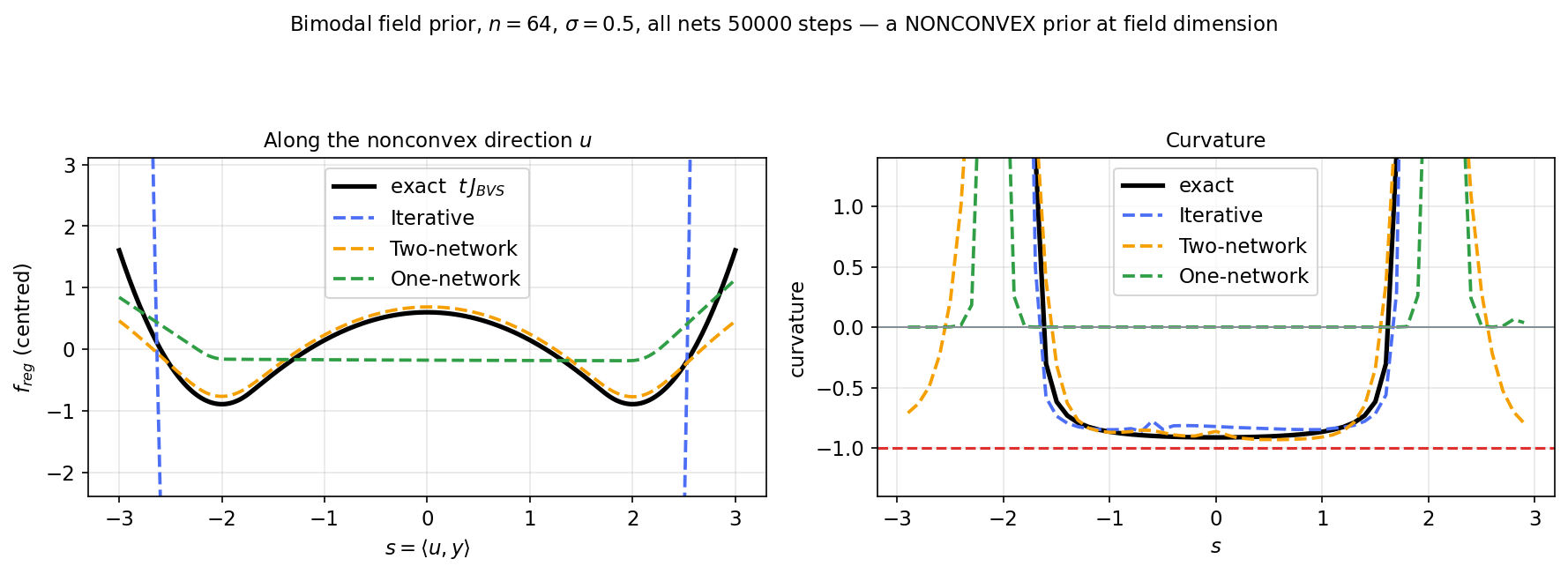}
 \caption{\Toadd{Bimodal field prior at $n=64$, along the nonconvex direction $s=\langle\bu,\by\rangle$. \emph{Left}: recovered regularizer, with the vertical axis clipped to the range of the exact curve. \emph{Right}: curvature. The one-network recovery is flat between the wells, and the iterative recovery leaves the operating support towards the ends of the interval.}}
 \label{fig:pnp-bimodal}
\end{figure}

In~\cref{tab:pnp-bimodal} the one-network recovery matches the $63$ Gaussian directions to $1.27\%$, eleven times better than the iterative recovery, and reads a curvature of $+0.0010$ along $\bu$ where the exact value is $-0.9113$. 

The iterative recovery attains a better value error, $1.56\%$ against $3.01\%$, and trails the two-network recovery by a factor of four on the gradient residual and by more than thirty on the perpendicular spectrum. Its inner solve converges, with certificate $\|\nabla\psiNN(\hat{\bx})-\by\|/\max(1,\|\by\|)$ of median $9.4\times10^{-4}$ and maximum $4.0\times10^{-3}$, yet it returns preimages reaching $\|\hat{\bx}\|_\infty=2.46$, outside the region where $\psiNN$ was fitted, and the recovered values diverge there. 

\Cref{fig:pnp-bimodal} shows the two semiconvex recoveries finding the dip. In
the left panel both track the exact $t\JBVS$ through the two wells and the local
maximum between them, while the one-network curve runs flat across the same
interval. The right panel repeats the curvature test of~\cref{fig:gaussian_mixture} at
field dimension, and the outcome is the same: both semiconvex curves reach the
plateau near $-0.91$, and the convex one does not leave zero.

We diagnose the two released operators~\cite{hu2023plug}, both of which fill the slot $\mathsf{D}$ in~\cref{eq:pirate-update}. The first, $\Dpir$, is PIRATE's operator, a convolutional network trained to remove Gaussian noise from registration fields. The second, $\Dplus$, is the same network after end-to-end fine-tuning of the whole registration pipeline. It does not itself denoise, since the fine-tuning targeted registration accuracy and was fit to no noise model. Neither is the ideal $\Dstar$, so~\cref{eq:tweedie} settles nothing about them and we measure the two Jacobian conditions directly.

At inputs drawn from the networks' own distribution we measure two quantities: the asymmetry ratio $\rho=\|\nabla \mathsf{D}-\nabla \mathsf{D}^{\top}\|_F/(2\|\nabla \mathsf{D}\|_F)$ and the extreme eigenvalues of the symmetric part $M=\tfrac12(\nabla \mathsf{D}+\nabla \mathsf{D}^{\top})$. We estimate both matrix-free, from Jacobian-vector and vector-Jacobian products. Run under identical settings on the trained convex-ICNN proximal operator of \cref{subsec:tv_exp}, whose architecture guarantees both conditions, the same diagnostics return $\rho=3.1\times10^{-13}$ and $\lambda(M)\in[0.082,1.000]$. A symmetric surrogate of the same convolutional architecture returns $\rho=0$, and the resolvable floor $\rho=1.1\times10^{-8}$ is set instead by the largest inconsistency between the Jacobian-vector and vector-Jacobian products on the operators themselves.

\begin{table}[!htbp]
 \centering
 \caption{\Toadd{Jacobian diagnosis of the two released operators, $8$ test points. A proximal map requires $\rho=0$ and $\lambda(M)\geqslant0$, and a convex regularizer additionally requires $\lambda(M)\leqslant1$. The resolvable asymmetry floor is $1.1\times10^{-8}$. Each $\rho$ is the mean over the eight points, constant across them to the digits shown.}}
 \label{tab:pnp-probe}
 \begin{tabular}{lrrrl}
 \hline
 operator & $\rho$ & $\lambda_{\min}(M)$ & $\lambda_{\max}(M)$ & verdict\\
 \hline
 ICNN prox & $3.1\times10^{-13}$ & $0.082$ & $1.000$ & proximal, convex\\
 $\Dpir$ & $0.44$ & $-0.72$ & $1.14$ & not a gradient field\\
 $\Dplus$ & $0.016$ & $0.95$ & $1.05$ & proximal, mildly nonconvex\\
 \hline
 \end{tabular}
\end{table}

In~\cref{tab:pnp-probe}, at $\rho=0.44$ no function has $\bx-\Dpir(\bx)$ as its gradient, so no $\freg$ has $\prox_{\freg}=\Dpir$, and the objective that PIRATE writes down does not exist for PIRATE's own denoiser. The test accepts $\Dplus$, at $\rho=0.016$, with symmetric part positive definite and spectrum in $[0.95,1.05]$.

The ratio $\rho$ normalizes by $\nabla \mathsf{D}$, while the regularizer's gradient is the residual $\bx-\mathsf{D}(\bx)$, whose Jacobian carries the same skew part over a smaller norm. For $\Dplus$ that residual is $2.2\%$ of $\nabla \Dplus$ in Frobenius norm, so measured in its own units the asymmetry is $\rho_{\mathrm{res}}=0.70$, and that skew part is an error no gradient field can remove. For $\Dpir$ the residual is $4.5$ times larger than $\nabla \Dpir$ and $\rho_{\mathrm{res}}=0.098$, which reverses the ordering, though its indefinite $M$ rules out a proximal reading either way.

We recover $\freg$ for $\Dplus$. On its operating distribution, the released registration field plus unit-variance Gaussian noise, $\Dplus$ is affine to $0.20\%$ over $12$ point pairs, with a Jacobian varying by $0.8\%$ between points. An affine proximal map has a quadratic regularizer, so writing $M$ for the symmetric part of the measured Jacobian and $\bm b$ for the constant,
\begin{equation}\label{eq:pirate-freg}
 \nabla\freg(\by)=M^{-1}(\by-\bm b)-\by,
\end{equation}
with no network trained. Its curvatures run from $-0.050$ to $+0.057$, read off the spectrum $[0.946,1.053]$ of $M$. At a held-out input the prox residual of~\cref{eq:pirate-freg} is $0.597$ against the ceiling $0.704$ that $\rho_{\mathrm{res}}$ and the affinity error impose, so the quadratic account is as accurate as the skew part permits. The same measurement on $\Dpir$ returns a pair error of $0.50$ and a Jacobian varying by $93\%$, so that operator is not affine and admits no such closed form.

The same two conditions form a reachability test. Given only that $\Dplus=\nabla_{\bx}\Keps(\cdot,t)$ with $\bx \mapsto \Keps(\bx,t)$ convex, the operator fixes a value function at its operating time $t=1$,
\begin{equation}\label{eq:pirate-S}
 \Seps(\bx,t)=\frac{1}{t}\left(\frac{1}{2}\normsq{\bx}-\Keps(\bx,t)\right),
 \qquad \Dplus=\mathrm{id}-t\,\nabla_{\bx}\Seps(\cdot,t),
\end{equation}
with no initial datum. Here $I_{t}(S)$ refers to the inviscid problem~\cref{eqn:intro2}, not to the viscous equation~\cref{eq:viscous_hj_tv} that produced $\Seps$: we read $\Seps(\cdot,t)$ as a candidate value function and ask which inviscid data reach it. Its Lipschitz hypothesis holds where we probe, since $t\nabla_{\bx}\Seps = \matr{I}_{n \times n} - \nabla\Dplus$ is bounded there by~\cref{tab:pnp-probe}. The criterion of~\cref{thm:reachability}, that $\bx\mapsto t\Seps(\bx,t)$ be $1$-semiconcave, then means $\nabla_{\bx}^{2}\Seps(\cdot,t)\preceq\matr{I}_{n \times n}/t$, which is the condition $\nabla \Dplus=\matr{I}_{n \times n}-t\nabla_{\bx}^{2}\Seps\succeq0$ that~\cref{tab:pnp-probe} measures. Here $\Dplus$ passes at $\lambda_{\min}=0.95>0$, so the forward flow reaches $\Seps(\cdot,t)$, while $\Dpir$ fails at $\lambda_{\min}=-0.72$ and is unreachable. By~\cref{thm:prior_optimality} the recovered $\JBVS=\freg/t$ is the initial datum that reaches $\Seps(\cdot,t)$, so $(\bx,t)\mapsto \Seps(\bx,t)$ is a Hamilton--Jacobi trajectory running from $\JBVS$ at time $0$ to the operator's own state at time $1$. The curvatures of~\cref{eq:pirate-freg} sit far above the floor $-1$ that $t=1$ imposes, so the flow clips nothing. 

\section{Discussion}\label{sec:discussion}
\GPLadd{We combined} the theory of viscosity solutions of HJ PDEs \GPLadd{with} recent work on proximal operators to develop deep learning \GPLadd{methods that learn the} prior of the proximal operator~\cref{eq:def_prox_argmin} from data. Our approach \GPLadd{builds} on connections between proximal operators and HJ PDEs, namely that proximal operators \GPLadd{arise} from solutions to a class of HJ PDEs, and on the theory \GPLadd{of} the inverse problem for HJ equations. As discussed in~\cref{sec:backward}, while \GPLadd{infinitely many priors may recover}~\cref{eq:def_prox_min}, \GPLadd{one choice is natural: reversing} the time in the HJ PDE~\cref{eq:finite_time_backward} and using the \GPLadd{value function} $(\bx,t) \mapsto S(\bx,t)$ as \GPLadd{terminal} condition. The resulting backward viscosity solution yields a semiconvex prior $\JBVS$ that can reconstruct the value function $(\bx,t) \mapsto S(\bx,t)$ underlying \GPLadd{the }proximal operator. \GPLadd{In}~\cref{sec:inverse_incomplete_info}\GPLadd{, we} considered the special case where only samples of the \GPLadd{proximal operator} and its values were available, and used max-plus algebra to derive characterizations and \GPLadd{error properties} of $\JBVS$ with respect to convex functions approximating it from \GPLadd{below}.

We also presented numerical experiments in~\cref{sec:numerics} to learn both the forward solution $\bx \mapsto S(\bx,t)$ and the prior $\JBVS$ by training a convex neural network\GPLadd{, specifically a learned proximal network,} from data $\{\bx_{k},S(\bx_{k},t)\}_{k=1}^{K}$\GPLadd{, and, where only the proximal operator is observed, from data $\{\bx_{k},\prox_{tJ}(\bx_{k})\}_{k=1}^{K}$}. \GPLadd{These numerical results} measure the relative $L^{2}$ accuracy of the recovered prior in dimensions up to $64$, where our method evaluates the prior in a single forward pass. \GPLadd{On the two families of~\cref{ex:quad_exact} a max-affine quadratic potential, whose conjugate is available in closed form, recovers the prior with no second network and no inversion. We also recovered the implicit regularizer of a total-variation posterior-mean denoiser, a function with no analytic form.} 
\GPLadd{Finally, we diagnose the two plug-and-play operators released in~\cite{hu2023plug}; the characterization shows that the residual of PIRATE's denoiser is not the gradient of a function. The fine-tuned operator satisfies both conditions discussed in~\cref{subsec:pnp_exp}, and its regularizer has a closed-form expression.}

While this work concerns proximal operators, we expect our approach can \GPLadd{extend} to a broad class of Bregman divergences, as recent results \GPLadd{on} inverse problems for HJ equations suggest~\cite{esteve2020inverse}. Another \GPLadd{direction would be to learn the prior} $J$ \GPLadd{from a known value function} $(\bx,t) \mapsto S(\bx,t)$ using Monte Carlo sampling strategies, as recently proposed in~\cite{park2025neural} for the forward problem of HJ equations (i.e., learning $(\bx,t) \mapsto S(\bx,t)$ from known $J$). It would also be interesting to devise similar deep learning methods for the inverse problem of HJ equations with possibly time- or state-dependent Hamiltonians, \GPLadd{a setting} relevant to optimal control problems.

\appendix
\section{Formal calculation of~\cref{eq:sol_hj_pde_pqm}}
\label{app:pqm_calc}
For each $k \in \{1,\dots,K\}$, define the quadratic
\[
q_{k}(\by) \coloneqq tJ(\by_{k}) + \frac{1}{2}\normsq{\bx_{k}-\by_{k}} - \frac{1}{2}\normsq{\bx_{k}-\by} + \frac{\GPLadd{\alpha'}}{2}\normsq{\by-\by_{k}},
\]
so that the explicit formula for $t\JPQM$ derived in \cref{subsec:piecewise_approxs} reads $t\JPQM(\by) = \max_{k \in \{1,\dots,K\}} q_{k}(\by)$, where we used $J(\by_{k}) = \JBVS(\by_{k})$ from \cref{thm:prior_optimality}. Fix $\bx \in \Rn$. Multiplying the objective of~\cref{eq:sol_hj_pde_pqm} by $t$ gives
\[
t\left(\frac{1}{2t}\normsq{\bx-\by} + \JPQM(\by)\right) = \max_{k \in \{1,\dots,K\}} F_{k}(\by), \qquad F_{k}(\by) \coloneqq \frac{1}{2}\normsq{\bx-\by} + q_{k}(\by).
\]
We first minimize each piece $F_{k}$. Expanding the squares gives
\[
\frac{1}{2}\normsq{\bx-\by} - \frac{1}{2}\normsq{\bx_{k}-\by} = \frac{1}{2}\normsq{\bx} - \frac{1}{2}\normsq{\bx_{k}} + \langle\bx_{k}-\bx,\by\rangle,
\]
and therefore
\[
F_{k}(\by) = tJ(\by_{k}) + \frac{1}{2}\normsq{\bx_{k}-\by_{k}} + \frac{1}{2}\normsq{\bx} - \frac{1}{2}\normsq{\bx_{k}} + \langle\bx_{k}-\bx,\by\rangle + \frac{\GPLadd{\alpha'}}{2}\normsq{\by-\by_{k}}.
\]
The function $F_{k}$ is $\GPLadd{\alpha'}$-strongly convex with $\nabla F_{k}(\by) = \bx_{k}-\bx + \GPLadd{\alpha'}(\by-\by_{k})$, so its unique minimizer is
\[
\hat{\by}_{k} = \by_{k} + \frac{\bx-\bx_{k}}{\GPLadd{\alpha'}}.
\]
Substituting $\hat{\by}_{k}$ into the two $\by$-dependent terms yields
\begin{align*}
\langle\bx_{k}-\bx,\hat{\by}_{k}\rangle + \frac{\GPLadd{\alpha'}}{2}\normsq{\hat{\by}_{k}-\by_{k}}
  &= \langle\bx_{k}-\bx,\by_{k}\rangle - \frac{1}{\GPLadd{\alpha'}}\normsq{\bx-\bx_{k}} + \frac{1}{2\GPLadd{\alpha'}}\normsq{\bx-\bx_{k}} \\
  &= \langle\bx_{k}-\bx,\by_{k}\rangle - \frac{1}{2\GPLadd{\alpha'}}\normsq{\bx-\bx_{k}},
\end{align*}
while expanding $\frac{1}{2}\normsq{\bx_{k}-\by_{k}}$ collapses the remaining terms:
\begin{align*}
\frac{1}{2}\normsq{\bx_{k}-\by_{k}} + \frac{1}{2}\normsq{\bx} - \frac{1}{2}\normsq{\bx_{k}} + \langle\bx_{k}-\bx,\by_{k}\rangle
  &= \frac{1}{2}\normsq{\by_{k}} + \frac{1}{2}\normsq{\bx} - \langle\bx,\by_{k}\rangle \\
  &= \frac{1}{2}\normsq{\bx-\by_{k}}.
\end{align*}
Combining the last three displays gives, for each $k \in \{1,\dots,K\}$,
\[
\min_{\by \in \Rn} F_{k}(\by) = tJ(\by_{k}) + \frac{1}{2}\normsq{\bx-\by_{k}} - \frac{1}{2\GPLadd{\alpha'}}\normsq{\bx-\bx_{k}}.
\]
It remains to pass from the pieces to their maximum. Since $\max_{j} F_{j} \geqslant F_{k}$ for every $k$, we have
\[
\inf_{\by \in \Rn}\left\{\frac{1}{2}\normsq{\bx-\by} + t\JPQM(\by)\right\} \geqslant \max_{k \in \{1,\dots,K\}} \min_{\by \in \Rn} F_{k}(\by).
\]
Suppose now that at a minimizer $\hat{\by}$ of the left-hand side a single piece is active, that is, $t\JPQM = q_{k^{*}}$ in a neighborhood of $\hat{\by}$ for some $k^{*}$. Then $\hat{\by}$ is a local, hence by convexity global, minimizer of $F_{k^{*}}$, so the left-hand side equals $\min_{\by} F_{k^{*}}(\by)$, and the inequality above collapses to an equality with the maximum attained at $k = k^{*}$. Dividing by $t$ yields~\cref{eq:sol_hj_pde_pqm} with this $k$. The single-active-piece assumption makes the calculation formal: at points $\bx$ whose minimizer $\hat{\by}$ lies on the boundary between two pieces, the left-hand side can exceed every $\min_{\by} F_{k}(\by)$, and~\cref{eq:sol_hj_pde_pqm} holds only as the lower bound displayed above.


\section*{Conflict of interest}
The authors declare that they have no conflict of interest.

\bibliographystyle{siamplain}
\bibliography{proj-bib}
\end{document}